\documentclass[reqno]{amsart}

\usepackage{amssymb, hyperref}

\begin{document}
\title[On well posedness for the INLS equations]{On well posedness for the inhomogeneous nonlinear Schr\"odinger equation}
	\author[C. M. GUZM\'AN ]
	{CARLOS M. GUZM\'AN }  
	
	\address{CARLOS M. GUZM\'AN \hfill\break
		Department of Mathematics, University Federal of Minas Gerais, BRAZIL}
	\email{carlos.guz.j@gmail.com}

\begin{abstract}
The purpose of this paper is to study well-posedness of the initial value problem (IVP) for the inhomogeneous nonlinear Schr\"odinger equation (INLS) 
$$
i u_t +\Delta u+\lambda|x|^{-b}|u|^\alpha u = 0, 
$$
where $\lambda=\pm 1$ and $\alpha$, $b>0$.\\
We obtain local and global results for initial data in $H^s(\mathbb{R}^N)$, with $0\leq s\leq 1$. To this end, we use the contraction mapping principle based on the Strichartz estimates related to the linear problem. 
		
\end{abstract}
	
	\maketitle  
	\numberwithin{equation}{section}
	\newtheorem{theorem}{Theorem}[section]
	\newtheorem{proposition}[theorem]{Proposition}
	\newtheorem{lemma}[theorem]{Lemma}
	\newtheorem{corollary}[theorem]{Corollary}
	\newtheorem{remark}[theorem]{Remark}
	\newtheorem{definition}[theorem]{Definition}

\section{Introduction}
\indent  In this work, we study the initial value problem (IVP), also called the Cauchy problem, for the inhomogenous nonlinear Schr\"odinger equation (INLS)
\begin{equation}\label{INLS}
\left\{\begin{array}{cl}
i\partial_tu +\Delta u + \lambda|x|^{-b} |u|^\alpha u =0, & \;\;\;t\in \mathbb{R} ,\;x\in \mathbb{R}^N,\\
u(0,x)=u_0(x), &
\end{array}\right.
\end{equation}
where $u = u(t,x)$ is a complex-valued function in space-time  $\mathbb{R}\times\mathbb{R}^N$, $\lambda=\pm 1$ and $\alpha, b>0$. The equation is called ``focusing INLS" when $\lambda= +1$ and ``defocusing INLS" when $\lambda= -1$. 
		
\ In the end of the last century, it was suggested that stable high power propagation can be achieved in a plasma by sending a preliminary laser beam that creates a channel with a reduced electron density, and thus reduces the nonlinearity inside the channel, see Gill \cite{GILL} and Liu-Tripathi \cite{LIU}. In this case, the beam propagation can be modeled by the inhomogeneous nonlinear Schr\"odinger equation in the following form:
$$
i\partial_tu +\Delta u + K(x) |u|^\alpha u =0.
$$
This model has been investigated by several authors, see, for instance, Merle \cite{MERLE} and Rapha\"el-Szeftel \cite{RAFAEL}, for $k_1<K(x)<k_2$ with $k_1,k_2>0$, and Fibich-Wang \cite{FIBI}, for $K(\epsilon |x|)$ with $\epsilon$ small and $K\in C^4(\mathbb{R}^N)\cap  L^{\infty}(\mathbb{R}^N)$. However, in these works $K(x)$ is bounded which is not verified in our case.
			
\ Our main goal here is to establish local and global results for the Cauchy problem \eqref{INLS} in $H^s(\mathbb{R}^N)$, with $0\leq s\leq 1$ applying Kato's method. Indeed, we construct a closed subspace of $C\left([-T,T];H^s(\mathbb{R}^N)\right)$  such that the operator defined by 
 \begin{equation}\label{OPERATOR} 
	G(u)(t)=U(t)u_0+i\lambda \int_0^t U(t-t')|x|^{-b}|u(t')|^\alpha u(t')dt',
\end{equation}
where $U(t)$ denotes the solution to the linear problem $i\partial_tu +\Delta u=0$, with initial data $u_0$, is stable and contractive in this space. Thus by the contraction mapping principle we obtain a unique fixed point. The fundamental tool to prove these results are the classic Strichartz estimates satisfied by the solution of the linear Schr\"odinger equation.

\ Notice that if $u(t,x)$ is solution of \eqref{INLS} so is $u_\delta(t,x)=\delta^{\frac{2-b}{\alpha}}u(\delta^2 t,\delta x)$, with initial data $u_{0,\delta}(x)$ for all $\delta >0$. Computing the homogeneus Sobolev norm we get
	$$
\|u_{0,\delta}\|_{\dot{H}^s}=\delta^{s-\frac{N}{2}+\frac{2-b}{\alpha}}\|u_0\|_{\dot{H}^s}.
$$
Hence, the scale-invariant Sobolev norm is $H^{s_c}(\mathbb{R}^N)$ with $s_c=\frac{N}{2}-\frac{2-b}{\alpha}$ (critical Sobolev index). Note that, if $s_c = 0$ (alternatively $\alpha = \frac{4-2b}{N}$) the problem is known as the mass-critical or $L^2$-critical; if $s_c=1$ (alternatively $\alpha =\frac{4-2b}{N-2}$) it is called energy-critical or $\dot{H}^1$-critical, finally the problem is known as mass-supercritical and energy-subcritical if $0<s_c<1$. On the other hand, the inhomogeneous nonlinear Schr\"odinger equation has the following conserved quantities:
\begin{equation}\label{mass}
Mass\equiv M[u(t)]=\int_{\mathbb{R}^N}|u(t,x)|^2dx=M[u_0]
\end{equation}
and
\begin{equation*}\label{energy}
Energy \equiv E[u(t)]=\frac{1}{2}\int_{\mathbb{R}^N}| \nabla u(t,x)|^2dx-\frac{\lambda}{\alpha +2}\left\| |x|^{-b}|u|^{\alpha +2}\right\|_{L^1_x}=E[u_0].
\end{equation*}

\ The well-posedness theory for the INLS equation \eqref{INLS} was studied for many authors in recent years. Let us briefly recall the best results available in the literature. Cazenave \cite{CAZENAVEBOOK} studied the well-posedness in $H^1(\mathbb{R}^N)$ using an abstract theory. To do this, he analyzed \eqref{INLS} in the sense of distributions, that is, $i\partial_tu +\Delta u + |x|^{-b} |u|^\alpha u =0$ in $H^{-1}(\mathbb{R}^N)$ for almost all $t\in I$. Therefore, using some results of Functional Analysis and Semigroups of Linear Operators, he proved  that it is appropriate to seek solutions of \eqref{INLS} satisfying      
$$
u \in  C \left([0, T);H^1(\mathbb{R}^N)\right) \cap C^1 \left([0, T);H^{-1}(\mathbb{R}^N)\right) \mbox { for some T} > 0.
$$
It was also proved that for the defocusing case ($\lambda=-1$) any local solution of the IVP (\ref{INLS}) with $u_0\in H^1(\mathbb{R}^N)$ extends globally in time.\\
 Other authors like Genoud-Stuart \cite{GENSTU} (see also references therein) also studied this problem for the focusing case ($\lambda=1$). Using the abstract theory developed by Cazenave \cite{CAZENAVEBOOK}, they showed that the IVP \eqref{INLS} is locally well-posed in $H^1(\mathbb{R}^N)$ if $0<\alpha <2^*$, where
\begin{equation}\label{def2*}
2^*:=\left\{\begin{array}{cl}
\frac{4-2b}{N-2}&\;\;\;\;N\geq 3,\\
\infty &\;\;\;\;\;N=1,2.
\end{array}\right.
\end{equation}
Recently, using some sharp Gagliardo-Nirenberg inequalities, Genoud \cite{GENOUD} and Farah \cite{LG} extended for the focusing INLS equation \eqref{INLS} some global well-posedness results obtained, respectively, by Weinstein \cite{WEINSTEIN} for the $L^2$-critical NLS equation and by Holmer-Roudenko \cite{HOLROU} for the $L^2$-supercritical and $H^1$-subcritical case. These authors proved that the solution $u$ of the Cauchy problem \eqref{INLS} is globally defined in $H^1(\mathbb{R}^N)$ quantifying the smallness condition in the initial data.

\ However, the abstract theory developed by Cazenave and later used by Genoud-Stuart \cite{GENSTU} to show well-posedness for \eqref{INLS}, does not give sufficient tools to study other interesting questions, for instance, scattering and  blow up investigated by Kenig-Merle \cite{KENIG}, Holmer-Roudenko-Duyckaerts \cite{DUCHOLROU} and others, for the NLS equation. To study these problems, the authors rely on the Strichartz estimates for NLS equation and the classical fixed point argument combining with the concentration-compactness and rigidity technique. 

\ Inspired by these papers and working toward the proof of scattering and blow up for the INLS equation, we show the well-posedness for the IVP \eqref{INLS} using the classic Strichartz estimates and the contraction mapping principle.

\ Applying this technique in the case $b=0$ (classical nonlinear Schr\"odinger equation (NLS)), the IVP \eqref{INLS} has been extensively studied over the three decades. The $L^2$-theory was obtained by Y. Tsutsumi \cite{TSUTSUMI} in the case $0<\alpha<\frac{4}{N}$. The $H^1$-subcritical case was studied by Ginibre-Velo \cite{GV}-\cite{GV1} and Kato \cite{KATO} (these papers also consider nonlinearities much more general than a pure power). Later, Cazenave-Weissler \cite{CAZEWE} treated the $L^2$-critical case and the $H^1$-critical case.

\ We summarize the well known well-posedness theory for the NLS equation in the following theorem (we refer, for instance, to Linares-Ponce \cite{FELGUS} for a proof of these results).  
\begin{theorem}\label{NLSMODEL} Consider the Cauchy problem for the NLS equation $($\eqref{INLS} with $b=0)$. Then, the following statements hold
\begin{enumerate}
\item  If $0<\alpha<\frac{4}{N}$, then the IVP \eqref{INLS} is locally and globally well posed in $L^2(\mathbb{R}^N)$. Moreover if $\alpha=\frac{4}{N}$, it is globally well posed in $L^2(\mathbb{R}^N)$ for small initial data.  
\item The IVP \eqref{INLS} with $b=0$ is locally well posed in $H^1(\mathbb{R}^N)$ if $0<\alpha \leq \frac{4}{N-2}$ for $N\geq 3$ or $0<\alpha <+\infty$, for $N=1,2$. Also, it is globally well-posed in $H^1(\mathbb{R}^N)$ if
\begin{itemize}
 \item [(i)] $\lambda<0$,
\item [(ii)] $\lambda>0$ and $ 0 < \alpha < \frac{4}{N}$,
 \item [(iii)] $\lambda>0$, $\frac{4}{N}<\alpha <\frac{4}{N-2}$ and small initial data,
\item [(iv)] $\lambda>0$, $\alpha =\frac{4}{N-2}$ and small initial data.
\end{itemize}
\end{enumerate}
\end{theorem}
\noindent In addition, Cazenave-Weissler \cite{CAZENAVEHS} and recently Cazenave-Fang-Han \cite{CAZENAVECONTINUOUS} showed that the IVP for the NLS is locally well posed in $H^s(\mathbb{R}^N)$ if $0 < \alpha\leq  \frac{4}{N-2s}$ and $0 < s < \frac{N}{2}$, moreover  the local solution extends globally in time for small initial data.

\ Our main interest in this paper is to prove similar results for the INLS equation. To this end, we divide in two parts.

\ The first part is devoted to study the local theory of the IVP \eqref{INLS}. We start considering the local well-posedness in $L^2(\mathbb{R}^N)$.


\begin{theorem}\label{LWPL2}
Let $0< \alpha<\frac{4-2b}{N}$ and $0<b<\min\{2,N\}$, then for all $u_0 \in L^2(\mathbb{R}^N)$ there exist $T=T(\|u_0\|_{L^2},N,\alpha)>0$ and a unique solution $u$ of \eqref{INLS} satisfying
$$
u \in C\left([-T,T];L^2(\mathbb{R}^N)\right) \cap L^q\left([-T,T];L^{r}(\mathbb{R}^N)\right), 
$$
for any ($q,r$) $L^2$-admissible. Moreover, the continuous dependence upon the initial data holds.
\end{theorem}
 It is worth to mention that the last theorem is an extension of a result by Tsutsumi \cite{TSUTSUMI} (which asserts local well-posedness for the NLS equation, \eqref{INLS} with $b=0$, when $0<\alpha <\frac{4}{N}$) to the INLS model.

\ Next, we treat the local well-posedness in $H^s(\mathbb{R}^N)$ for $0<s\leq 1$. Before stating the theorem, we define the following numbers
\begin{equation}\label{def2s}
\widetilde{2}:= \left\{\begin{array}{cl}
\frac{N}{3}&\;\;N=1,2,3,\\
2&\;\;N\geq 4
\end{array}\right.\;\;\;\textnormal{and}\;\;\;\;  \alpha_s:=\left\{\begin{array}{cl}
\frac{4-2b}{N-2s}&\;\;s<\frac{N}{2},\\
+\infty & \;\;s=\frac{N}{2}.
\end{array}\right.
\end{equation}
\begin{theorem}\label{LWPHs} 
Assume $ 0<\alpha<\alpha_s$, $0<b<\widetilde{2}$ and $\max\{0,s_c\}<s\leq \min\{\frac{N}{2},1\}$. If $u_0 \in H^s(\mathbb{R}^N)$ then there exist $T=T(\|u_0\|_{H^s},N,\alpha)>0$ and a unique solution $u$ of \eqref{INLS} with
$$
u \in C\left([-T,T];H^s(\mathbb{R}^N) \right) \cap L^q\left([-T,T];H^{s,r}(\mathbb{R}^N)    \right) 
$$
for any ($q,r$) $L^2$-admissible. Moreover, the continuous dependence upon the initial data holds.
\end{theorem}
\begin{remark}
 Observe that $\alpha<\frac{4-2b}{N-2s}$ is equivalent to $s_c<s$. On the other hand, if $0<\alpha<\frac{4-2b}{N}$ then $s_c<0$, for this reason we add the restriction $s>\max\{0,s_c\}$ in the above statement. 
\end{remark}

\ As an immediate consequence of the Theorem \ref{LWPHs}, we have that the IVP \eqref{INLS} is locally well-posed in $H^1(\mathbb{R}^N)$. 
\begin{corollary} \label{LWPH1}
Assume $N\geq 2$, $0<\alpha<2^*$ and $0<b<\widetilde{2}$. If $u_0 \in H^1(\mathbb{R}^N)$ then the initial value problem \eqref{INLS} is locally well-posed and
$$
u \in C\left([-T,T];H^1(\mathbb{R}^N) \right) \cap L^q\left([-T,T];H^{1,r}(\mathbb{R}^N)    \right),
$$
for any ($q,r$) $L^2$-admissible.
\end{corollary}
\begin{remark} One important difference of the previous results and its its counterpart for the NLS model (see Theorem \ref{NLSMODEL}-(2)) is that we do not treat the critical case here, i.e. $\alpha=\frac{4-2b}{N-2s}$ with $0\leq s\leq 1$ and $N\geq 3$. It is still an open problem. 
\end{remark}

\ In the second part, we consider the global well-posedness of the Cauchy problem \eqref{INLS}. We begin with a global result in $L^2(\mathbb{R}^N)$ which is an immediate consequence of Theorem \ref{LWPL2}.  

\begin{theorem}\label{GWPL2}
If $0< \alpha<\frac{4-2b}{N}$ and $0<b<\min\{2,N\}$, then for all $u_0 \in L^2(\mathbb{R}^N)$ the local solution $u$ of the IVP \eqref{INLS} extends globally with
$$
u \in C\left(\mathbb{R};L^2(\mathbb{R}^N)\right) \cap L^q\left(\mathbb{R};L^{r}(\mathbb{R}^N)\right), 
$$
for any ($q,r$) $L^2$-admissible. 
\end{theorem}

\ In the sequel we establish a small data global theory for the INLS model \eqref{INLS}.

\begin{theorem}\label{GWPHs} 
Let $ \frac{4-2b}{N}<\alpha<\alpha_s$ with $0<b<\widetilde{2}$, $s_c<s\leq \min\{\frac{N}{2},1\}$ and $u_0 \in H^s(\mathbb{R}^N)$. If $\|u_0\|_{H^s}\leq A$ then there exists $\delta=\delta(A)$ such that if $\|U(t)u_0\|_{S(\dot{H}^{s_c})}<\delta$, then the solution of \eqref{INLS} is globally defined. Moreover,
\begin{equation*}\label{NGWP1} 
\|u\|_{S(\dot{H}^{s_c})}\leq 2\|U(t)u_0\|_{S(\dot{H}^{s_c})} 
\end{equation*}
and
\begin{equation*}\label{NGWP2}
\|u\|_{S\left(L^2\right)}+\|D^s  u\|_{S\left(L^2\right)}\leq 2c\|u_0\|_{H^s}.
\end{equation*}
\end{theorem}
\begin{remark}
Note that in the last result we don't need the condition $s>\max\{0,s_c\}$ as in Theorem \ref{LWPHs}, since $\alpha>\frac{4-2b}{N}$ implies $s_c>0$.
\end{remark} 
\begin{remark}
Also note that by the Strichartz estimates \eqref{SE2}, the condition $\|U(t)u_0\|_{S(\dot{H}^{s_c})}<\delta$ is automatically satisfied if $\|u_0\|_{\dot{H}^{s_c}}\leq \frac{\delta}{c}$.
\end{remark}

\ A similar small data global theory for the NLS model can be found in Cazenave-Weissler \cite{CAZWEI}, Holmer-Roudenko \cite{HOLROU} and Guevara \cite{GUEVARA}. A consequence of the Theorem \ref{GWPHs} is the following global well-posed result in $H^1(\mathbb{R}^N)$.       
			
\begin{corollary}\label{GWPH1}
Let $N\geq 2$, $\frac{4-2b}{N}<\alpha<2^*$ with $0<b<\widetilde{2}$ and $u_0 \in H^1(\mathbb{R}^N)$. Assume $\|u_0\|_{H^1}\leq A$ then there exists $\delta=\delta(A)>0$ such that if $\|U(t)u_0\|_{S(\dot{H}^{s_c})}<\delta$, then  there exists a unique global solution $u$ of \eqref{INLS} such that
\begin{equation*}\label{NGWP3}
\|u\|_{S(\dot{H}^{s_c})}\leq  2\|U(t)u_0\|_{S(\dot{H}^{s_c})} 
\end{equation*}
and
\begin{equation*}\label{NGWP4}
\|u\|_{S\left(L^2\right)}+\|\nabla  u\|_{S\left(L^2\right)}\leq 2c\|u_0\|_{H^1}.
\end{equation*}
\end{corollary}
			
\ The rest of the paper is organized as follows. In section $2$, we introduce some notations and give a review of the Strichartz estimates. In section $3$, we prove the local well-posedness results: Theorems \ref{LWPL2} and \ref{LWPHs}. Finally, in Section $4$, we prove the results concerning the global theory: Theorems \ref{GWPL2} and \ref{GWPHs}.

\section{Notation and preliminares}
			
\ Let us start this section by introducing the notation used throughout the paper. We use $c$ to denote various constants that may vary line by line. Let a set $A\subset \mathbb{R}^N$, $A^C=\mathbb{R}^N \backslash A$ denotes the complement of $A$. Given $x,y \in \mathbb{R}^N$, $x.y$ denotes the inner product of $x$ and $y$ on $\mathbb{R}^N$. 

\ Let $ q,r \geq 1$, $T > 0$ and $s\in \mathbb{R}$, the mixed norms in the spaces $L^q_{[0,T]}L^r_x$ and $L^q_{[0,T]} H^s_x$ of $f(x,t)$ are defined, respectively, as
$$
\|f\|_{L^q_{0,T}L^r_x}=\left(\int_0^T\|f(t,.)\|^q_{L^r_x}dt\right)^{\frac{1}{q}}
$$
and
$$
\|f\|_{L^q_{0,T}H^s_x}=\left(\int_0^T\|f(t,.)\|^q_{H^s_x}dt\right)^{\frac{1}{q}}
$$
with the usual modifications when\footnote{$\|f\|_{L^\infty_{0,T}}=\sup\limits_{t\in [0,T]}|f(t)|$.} $q=\infty$ or $r=\infty$. In the case when $I=[0,T]$ and we restrict the $x$-integration to a subset $A\subset\mathbb{R}^N$ then the mixed norm will be denoted by $\|f\|_{L_I^qL^r_x(A)}$. Moreover, when $f(t, x)$ is defined for every time $t \in \mathbb{R}$ we shall consider the notations $\|f\|_{L_t^qL^r_x}$ and $\|f\|_{L_t^qH^s_x}$.

\ For $s\in \mathbb{R}$, $J^s$ and $D^s$ denote the Bessel and the Riesz potentials of order $s$, given via Fourier transform by the formulas
$$
 \widehat{J^s f}=(1+|\xi|^2)^{\frac{s}{2}}\widehat{f}\;\;\;\textnormal{and} \;\;\;\;\widehat{D^sf}=|\xi|^s\widehat{f},
 $$
where the Fourier transform of $f(x)$ is given by
$$
\widehat{f}(y)=\int_{\mathbb{R}^N}e^{-ix.\xi}f(x)dx.
$$
On the other hand, we define the norm of the Sobolev spaces $H^{s,r}(\mathbb{R}^N)$ and $\dot{H}^{s,r}(\mathbb{R}^N)$, respectively, by
$$
\|f\|_{H^{s,r}}:=\|J^sf\|_{L^r}\;\;\;\;\textnormal{and}\;\;\;\;\|f\|_{\dot{H}^{s,r}}:=\|D^sf\|_{L^r}.
$$
If $r=2$ we denote $H^{s,2}$ simply by $H^s$.

\ Next, we recall some Strichartz type estimates associated to the linear Schr\"odinger propagator.\\
\textbf{Strichartz type estimates.} We say the pair $(q,r)$ is $L^2$-admissible or simply admissible par if they satisfy the condition\footnote{We included in the above definition the improvement, due
to M. Keel and T. Tao \cite{tao:keel}, to the limiting case for Strichartz’s inequalities.}
\begin{equation*}
\frac{2}{q}=\frac{N}{2}-\frac{N}{r},
\end{equation*}
where 
\begin{equation}\label{L2Admissivel}
\left\{\begin{array}{cl}
2\leq & r  \leq \frac{2N}{N-2}\hspace{0.5cm}\textnormal{if}\;\;\;  N\geq 3,\\
2 \leq  & r < +\infty\;  \hspace{0.5cm}\textnormal{if}\;\; \;N=2,\\
2 \leq & r \leq + \infty\;  \hspace{0.5cm}\textnormal{if}\;\;\;N=1.
\end{array}\right.
\end{equation}
We also called the pair $\dot{H}^s$-admissible if
\begin{equation}\label{CPA1}
\frac{2}{q}=\frac{N}{2}-\frac{N}{r}-s,
\end{equation}
where 
\begin{equation}\label{CPA2}
\left\{\begin{array}{cl}
\frac{2N}{N-2s}\leq & r  \leq\left(\frac{2N}{N-2}\right)^{-}\;\hspace{0.4cm}\textnormal{if}\;\;\;  N\geq 3,\\
\frac{2}{1-s} \leq  & r \leq \left((\frac{2}{1-s})^+\right)'\;  \hspace{0.2cm}\textnormal{if}\;\; \;N=2,\\
\frac{2}{1-2s} \leq & r \leq + \infty\; \; \hspace{1.2cm}\textnormal{if}\;\;\;N=1.
\end{array}\right.
\end{equation}
Here, $a^-$ is a fixed number slightly smaller than a ($a^-=a-\varepsilon$ with $\varepsilon>0$ small enough) and, in a similar way, we define $a^+$. Moreover $(a^+)'$ is the number such that 
\begin{equation}\label{a^+}
\frac{1}{a}=\frac{1}{(a^+)'}+\frac{1}{a^+},
\end{equation}
that is $(a^+)':=\frac{a^+.a}{a^+-a}$ with $a^+$. Finally we say that $(q,r)$ is $\dot{H}^{-s}$-admissible if 
$$
\frac{2}{q}=\frac{N}{2}-\frac{N}{r}+s,
$$
where
\begin{equation}\label{H-s}
\left\{\begin{array}{cl}
\left(\frac{2N}{N-2s}\right)^{+}\leq & r  \leq\left(\frac{2N}{N-2}\right)^{-}\;\;\hspace{0.4cm}\textnormal{if}\;\;  N\geq 3,\\
\left(\frac{2}{1-s}\right)^{+} \leq  & r \leq \left((\frac{2}{1+s})^+\right)'\;  \hspace{0.2cm}\textnormal{if}\;\; \;N=2,\\
\left(\frac{2}{1-2s}\right)^{+} \leq & r \leq + \infty\; \; \hspace{1.2cm}\textnormal{if}\;\;\;N=1.
\end{array}\right.
\end{equation}

\  Given $s\in \mathbb{R}$, let $\mathcal{A}_s=\{(q,r);\; (q,r)\; \textnormal{is} \;\dot{H}^s-\textnormal{admissible}\}$ and $(q',r')$ is such that $\frac{1}{q}+\frac{1}{q'}=1$ and $\frac{1}{r}+\frac{1}{r'}=1$ for $(q,r)\in \mathcal{A}_s$. We define the following Strichartz norm
$$
\|u\|_{S(\dot{H}^{s})}=\sup_{(q,r)\in \mathcal{A}_{s}}\|u\|_{L^q_tL^r_x} 
$$
and the dual Strichartz norm
$$
\|u\|_{S'(\dot{H}^{-s})}=\inf_{(q,r)\in \mathcal{A}_{-s}}\|u\|_{L^{q'}_tL^{r'}_x}.
$$
Note that, if $s=0$ then $\mathcal{A}_0$ is the set of all $L^2$-admissible pairs. Moreover, if $s=0$, $S(\dot{H}^0)=S(L^2)$ and $S'(\dot{H}^{0})=S'(L^2)$. We just write $S(\dot{H}^s)$ or $S'(\dot{H}^{-s})$ if the mixed norm is evaluated over $\mathbb{R}\times\mathbb{R}^N$. To indicate a restriction to a time interval $I\subset (-\infty,\infty)$ and a subset $A$ of $\mathbb{R}^N$, we will consider the notations $S(\dot{H}^s(A);I)$ and $S'(\dot{H}^{-s}(A);I)$. 

\ We now list (without proving) some estimates that will be useful in our work.  
 \begin{lemma}\textbf{(Sobolev embedding)}\label{SI} Let $s\in (0,+\infty)$ and $1\leq p<+\infty$.
\begin{itemize}
\item [(i)] If $s\in (0,\frac{N}{p})$ then $H^{s,p}(\mathbb{R}^N)$ is continuously embedded in $L^r(\mathbb{R^N})$ where $s=\frac{N}{p}-\frac{N}{r}$. Moreover, 
\begin{equation}\label{SEI} 
\|f\|_{L^r}\leq c\|D^sf\|_{L^{p}}.
\end{equation}
\item [(ii)] If $s=\frac{N}{2}$ then $H^{s}(\mathbb{R}^N)\subset L^r(\mathbb{R^N})$ for all $r\in[2,+\infty)$. Furthermore,
\begin{equation}\label{SEI1} 
\|f\|_{L^r}\leq c\|f\|_{H^{s}}.
\end{equation}
\end{itemize}
\begin{proof} See Bergh-L\"ofstr\"om \cite[Theorem $6.5.1$]{BERLOF} (see also Linares-Ponce \cite[Theorem $3.3$]{FELGUS} and Demenguel-Demenguel \cite[Proposition 4.18]{DEMENGEL}). 
\end{proof}
\end{lemma}
\begin{remark}\label{SEI2} Using $(i)$, with $p=2$, we have that $H^s(\mathbb{R}^N)$, with $s\in (0,\frac{N}{2})$, is continuously embedded in $L^r(\mathbb{R}^N)$ and 
\begin{equation}\label{SEI22} 
\|f\|_{L^r}\leq c\|f\|_{H^{s}},
\end{equation}    
where $r\in[2,\frac{2N}{N-2s}]$.
\end{remark}
\begin{lemma}\textbf{(Fractional product rule)}\label{PRODUCTRULE}
Let $s\in (0, 1]$ and $1 < r, r_1, r_2, p_1, p_2 <+\infty$ are such that $\frac{1}{r}=\frac{1}{r_i}+\frac{1}{p_i}$ for $i=1,2$. Then,
$$
\| D^s (fg) \|_{L^r} \leq c \|f\|_{L^{r_1}}\|D^s g\|_{L^{p_1}}+c \|D^sf\|_{L^{r_2}}\| g\|_{L^{p_2}}.
$$
\begin{proof} See Christ-Weinstein \cite[Proposition $3.3$]{CHRWEI}.
\end{proof}
\end{lemma}
\begin{lemma}\textbf{(Fractional chain rule)}\label{CHAINRULE}
Suppose $G\in C^1(\mathbb{C})$, $s\in (0, 1]$, and $1 < r, r_1, r_2 <+\infty$ are such that $\frac{1}{r}=\frac{1}{r_1}+\frac{1}{r_2}$. Then,
$$
\| D^s G(u) \|_{L^r} \leq c \|G^{'}(u)\|_{L^{r_1}}\|D^s u\|_{L^{r_2}}.
$$
\begin{proof} See Christ-Weinstein \cite[Proposition $3.1$]{CHRWEI}.
\end{proof}
\end{lemma}
				
 \ The main tool to show the local and global well-posedness are the well-known Strichartz estimates. See for instance Linares-Ponce \cite{FELGUS} and Kato \cite{KATO} (see also Holmer-Roudenko \cite{HOLROU} and Guevara \cite{GUEVARA}).
\begin{lemma}\label{Lemma-Str} The following statements hold.
 \begin{itemize}
\item [(i)] (Linear estimates).
\begin{equation}\label{SE1}
\| U(t)f \|_{S(L^2)} \leq c\|f\|_{L^2},
\end{equation}
\begin{equation}\label{SE2}
\|  U(t)f \|_{S(\dot{H}^s)} \leq c\|f\|_{\dot{H}^s}.
\end{equation}
\item[(ii)] (Inhomogeneous estimates).
\begin{equation}\label{SE3}					 
\left \| \int_{\mathbb{R}} U(t-t')g(.,t') dt' \right\|_{S(L^2)}\;+\; \left \| \int_{0}^t U(t-t')g(.,t') dt' \right \|_{S(L^2) } \leq c\|g\|_{S'(L^2)},
\end{equation}
\begin{equation}\label{SE5}
\left \| \int_{0}^t U(t-t')g(.,t') dt' \right \|_{S(\dot{H}^s) } \leq c\|g\|_{S'(\dot{H}^{-s})}.
\end{equation}
\end{itemize}
\end{lemma} 
The relations \eqref{SE3} and \eqref{SE5} will be very useful to perform estimates on the nonlinearity $|x|^{-b}|u|^\alpha u$.

\ We end this section with three important remarks. 	
 \begin{remark}\label{nonlinerity}
Let $F(x,u)=|x|^{-b}|z|^\alpha z$, where $f(z)=|z|^\alpha z$. The complex derivative of $f$ is
\begin{equation*}
f_z(z)=\frac{\alpha+2}{2}|z|^\alpha\;\;\;\;\;\textnormal{and}\;\;\;\; f_{\bar{z}}    (u)=\frac{\alpha}{2}|z|^{\alpha-2}z^2. 
\end{equation*}
For $z,w\in \mathbb{C}$, we have 
\begin{equation*}
 f(z)-f(w)=\int_{0}^{1}\left[f_z(w+t(z-w))(z-w)+f_{\bar{z}}(w+t(z-w))\overline{(z-w)}\right]dt.
\end{equation*}
Thus,
\begin{equation}\label{FEI}
 |F(x,z)-F(x,w)|\lesssim |x|^{-b}\left( |z|^\alpha+ |w|^\alpha \right)|z-w|.
\end{equation}
\end{remark}

\begin{remark}\label{RB} 
Let $B=B(0,1)=\{ x\in \mathbb{R}^N;|x|\leq 1\}$ and $b>0$. If $x\in B^C$ then $|x|^{-b}<1$ and so
$$ 
\left	\||x|^{-b}f \right\|_{L^r_x}\leq \|f\|_{L_x^r(B^C)}+\left\||x|^{-b}f\right\|_{L_x^r(B)}. 
$$
\end{remark}

\ The next remark provides a condition for the integrability of $|x|^{-b}$ on $B$ and $B^C$. 
\begin{remark}\label{RIxb} 
We notice that if $\frac{N}{\gamma}-b>0$ then $\||x|^{-b}\|_{L^\gamma(B)}<+\infty$, indeed 
\begin{equation*}
\int_{B}|x|^{-\gamma b}dx=c\int_{0}^{1}r^{-\gamma b}r^{N-1}dr=c_1 \left. r^{N-\gamma b} \right |_0^1<+\infty\;\;\textnormal{if}\;\;N-\gamma b>0.
\end{equation*}
Similarly, we have that $\||x|^{-b}\|_{L^\gamma(B^C)}$ is finite if $\frac{N}{\gamma}-b<0$.
\end{remark}

\section{Local well-posedness}
				
\ In this section we prove the local well-posedness results. The theorems follows from a contraction mapping argument based on the Strichartz estimates. First, we show the local well-posedness in $L^2(\mathbb{R}^N)$ (Theorem \ref{LWPL2}) and then in $H^s(\mathbb{R}^N)$ for $0<s\leq 1$ (Theorem \ref{LWPHs}) as well as Corollary \ref{LWPH1}. 

\subsection{$L^2$-Theory} 
We begin with the following lemma. It provides an estimate for the INLS model nonlinearity in the Strichartz spaces. 
\begin{lemma}\label{lemmaL2} 
Let $0<\alpha<\frac{4-2b}{N}$  and $0<b<\min \{2,N\}$. Then,
\begin{equation}\label{NlemmaL2} 
\left\||x|^{-b}|u|^{\alpha}v \right\|_{S'(L^2;I)}\leq c (T^{\theta_1}+T^{\theta_2})\| u\|^{\alpha}_{S(L^2;I)}\| v\|_{S(L^2;I)},
\end{equation}
where $I=[0,T]$ and $c,\theta_1,\theta_2 >0$.
\begin{proof}
\ By Remark \ref{RB}, we have
\begin{eqnarray*}
\left\||x|^{-b}|u|^{\alpha}v \right\|_{S'\left(L^2;I\right)} &\leq & \left\||u|^{\alpha}v
\right\|_{S'\left(L^2(B^C);I\right)}+ \left\||x|^{-b}|u|^{\alpha}v \right\|_{S'\left(L^2(B);I\right)}\\
&\equiv& A_1+A_2. 
\end{eqnarray*}
Note that in the norm $A_1$ we don't have any singularity, so we know that 
\begin{equation}\label{A1} %
A_1\leq  cT^{\theta_1}\| u\|^{\alpha}_{S(L^2;I)}\| v\|_{S(L^2;I)}, 	
\end{equation}
where $\theta_1>0$. See Kato \cite[Theorem $0$]{KATO} (also see Linares-Ponce \cite[Theorem $5.2$ and Corollary $5.1$]{FELGUS}).
				
\ On the other hand, we need to find an admissible pair to estimate $A_2$. In fact, using the H\"older inequality twice we obtain 
\begin{eqnarray*}
A_2 &\leq & \left\|  |x|^{-b}|u|^{\alpha} v \right\|_{L^{q'}_IL_x^{r'}(B)} \leq \left \| \||x|^{-b}\|_{L^\gamma(B) } \|u\|^{\alpha}_{L_x^{\alpha r_1}} \|v\|_{L_x^{ r}}\right\|_{L_I^{q'}}\\
&\leq&  \||x|^{-b}\|_{L^\gamma(B)}T^{\frac{1}{q_1}} \| u\|^{\alpha}_{L_I^{\alpha q_2}L_x^{\alpha r_1}}\| v\|_{L_I^{q}L_x^{r}}\\
&\leq&  T^{\frac{1}{q_1}}\||x|^{-b}\|_{L^\gamma(B)} \| u\|^{\alpha}_{L_I^{q}L_x^r}\| v\|_{L_I^{q}L_x^{r}},
\end{eqnarray*}			
if $(q,r)$ $L^2$-admissible and 
\begin{equation}\label{lemmaL21}
\left\{\begin{array}{cl}\vspace{0.1cm}
\frac{1}{r'}=&\frac{1}{\gamma}+\frac{1}{r_1}+\frac{1}{r} \\ \vspace{0.1cm}
\frac{1}{q'}=&\frac{1}{q_1}+\frac{1}{q_2}+\frac{1}{q}\\
q=& \alpha q_2\;\;,\;\;r=\alpha r_1.
\end{array}\right.
\end{equation}
In order to have $\||x|^{-b}\|_{L^\gamma(B)}<+\infty$ we need $\frac{N}{\gamma}>b$, by Remark \ref{RIxb}. Hence, in view of \eqref{lemmaL21} $(q,r)$ must satisfy
\begin{equation}\label{lemmaL22}
\left\{\begin{array}{cl}\vspace{0.1cm}
 \frac{N}{\gamma}=&N-\frac{N(\alpha+2)}{r}>b\\
\frac{1}{q_1}=&1-\frac{\alpha+2}{q}.
\end{array}\right.
 \end{equation}
From the first equation in \eqref{lemmaL22} we have $N-b-\frac{N(\alpha+2)}{r}>0$, which is equivalent to 
\begin{equation}\label{lemmaL23}
\alpha<\frac{r(N-b)-2N}{N},
\end{equation}
for $r>\frac{2N}{N-b}$. By hypothesis $\alpha <\frac{4-2b}{N}$, then setting $r$ such that 
$$
\frac{r(N-b)-2N}{N}=\frac{4-2b}{N},
$$ 
we get\footnote{Since $0<b<\min\{N,2\}$ the denominator of $r$ is positive and $r>\frac{2N}{N-b}$. Moreover, by a simple computations we have $2\leq r\leq \frac{2N}{N-2}$ if $N\geq 3$, and $2\leq r <+\infty$ if $N=1,2$, that is $r$ satisfies \eqref{L2Admissivel}. Therefore, the pair $(q,r)$ above defined is $L^2$-admissible.} $r=\frac{4-2b+2N}{N-b}$ satisfying \eqref{lemmaL23}. Consequently, since $(q,r)$ is $L^2$-admissible we obtain $q=\frac{4-2b+2N}{N}$. Next, applying the second equation in \eqref{lemmaL22} we deduce 
$$
\frac{1}{q_1}=\frac{4-2b-\alpha N}{4-2b+2N},
$$ 
which is positive by the hypothesis $\alpha<\frac{4-2b}{N}$. Thus, 
\begin{equation*}\label{A2}
A_2\leq c T^{\theta_2} \| u\|^{\alpha}_{S(L^2;I)}\| v\|_{S(L^2;I)},
\end{equation*}
where $\theta_2=\frac{1}{q_1}$.
Therefore, combining \eqref{A1} and the last inequality we prove \eqref{NlemmaL2}.	
\end{proof}
\end{lemma}

\ Our goal now is to show Theorem \ref{LWPL2}.
					
\begin{proof}[\bf{Proof of Theorem \ref{LWPL2}}]  
We define $$X= C\left( [-T,T];L^2(\mathbb{R}^N)\right) \bigcap  L^q\left([-T,T];L^{r}(\mathbb{R}^N)\right),$$ for any ($q,r$) $L^2$-admissible, and 
\begin{equation*}\label{BL} 
 B(a,T)=\{u \in X : \|u\|_{S\left(L^2;[-T,T]\right)}\leq a \},
\end{equation*}
where $a$ and $T$ are positive constants to be determined later. We follow the standard fixed point argument to prove this result. It means that for appropriate values of $a$, $T$ we shall show that $G$ defined in \eqref{OPERATOR} defines a contraction map on $B(a,T)$.

\ Without loss of generality we consider only the case $t>0$. Applying Strichartz inequalities (\ref{SE1}) and (\ref{SE3}), we have
\begin{equation}\label{NSDL2} 
\|G(u)\|_{S(L^2;I)}\leq c\|u_0\|_{L^2}+ c\||x|^{-b}|u|^{\alpha+1}  \|_{S'(L^2;I)},
\end{equation}
where $I=[0,T]$. Moreover, Lemma \ref{lemmaL2} yields
\begin{eqnarray*}
\|G(u)\|_{S(L^2;I)}&\leq &c\|u_0\|_{L^2}+c (T^{\theta_1}+T^{\theta_2})\| u\|^{\alpha+1}_{S(L^2;I)}\\
&\leq& c\|u_0\|_{L^2}+c (T^{\theta_1}+T^{\theta_2}) a^{\alpha+1},
\end{eqnarray*}
provided $u\in B(a,T)$. Hence,
\begin{eqnarray*}
\|G(u)\|_{S\left(L^2;[-T,T]\right)} &\leq& c\|u_0\|_{L^2}+c (T^{\theta_1}+T^{\theta_2}) a^{\alpha+1}.
\end{eqnarray*}
Next, choosing $a=2c\|u_0\|_{L^2}$ and $T>0$ such that 
\begin{equation}\label{NT3}
c a^{\alpha} (T^{\theta_1}+T^{\theta_2}) < \frac{1}{4},
\end{equation}
we conclude  $G(u)\in B(a,T)$. 

\ Now we prove that $G$ is a contraction. Again using Strichartz inequality (\ref{SE3}) and \eqref{FEI}, we deduce
\begin{eqnarray*}
\|G(u)-G(v)\|_{S(L^2;I)}&\leq& c \left\||x|^{-b}(|u|^{\alpha}u-|v|^\alpha v)\right\|_{S'(L^2;I)}\\
&\leq& c\left\| |x|^{-b} |u|^{\alpha}|u-v|\right\|_{S'(L^2;I)}\\
& &+\;c\left\| |x|^{-b} |v|^{\alpha}|u-v|\right\|_{S'(L^2;I)}\\
&\leq& c (T^{\theta_1}+T^{\theta_2})\|u\|^\alpha_{S(L^2;I)}\|u-v\|_{S(L^2;I)}\\
&& +c (T^{\theta_1}+T^{\theta_2})\|v\|^\alpha_{S(L^2;I)} \|u-v\|_{S(L^2;I)},
\end{eqnarray*}
where $I=[0,T]$. That is,
\begin{align*}
\|G(u)-G(v)\|_{S(L^2;I)}&\leq c (T^{\theta_1}+T^{\theta_2})\left(\|u\|^\alpha_{S(L^2;I)} +\|v\|^\alpha_{S(L^2;I)} \right)\|u-v\|_{S(L^2;I)}\\
& \leq 2c (T^{\theta_1}+T^{\theta_2})a^\alpha \|u-v\|_{S(L^2;I)},
\end{align*}
provided $u,v\in B(a,T)$. Therefore, the inequality $\eqref{NT3}$ implies that
\begin{eqnarray*}
\|G(u)-G(v)\|_{S\left(L^2;[-T,T]\right)} &\leq & 2c (T^{\theta_1}+T^{\theta_2})a^\alpha \|u-v\|_{S\left(L^2;[-T,T]\right)}\\
&< & \frac{1}{2}\|u-v\|_{S\left(L^2;[-T,T]\right)},
\end{eqnarray*}
i.e., $G$ is a contraction on $S(a,T)$.
					
\ The proof of the continuous dependence is similar to the one given above and it will be omitted.
\end{proof}

\subsection{$H^s$-Theory}
The aim of this subsection is to prove the local well-posedness in $H^s(\mathbb{R}^N)$ with $0<s\leq 1$ (Theorem \ref{LWPHs}) as well as Corollary \ref{LWPH1}. Before doing that we establish useful estimates for the nonlinearity $|x|^{-b}|u|^\alpha u$. First, we consider the nonlinearity in the space $S'(L^2)$ and in the sequel in the space $D^{-s}S'(L^2)$, that is, we estimate the norm $\left\||x|^{-b}|u|^\alpha u \right\|_{S'(L^2;I)}$ and $\left\| D^s(|x|^{-b}|u|^\alpha u )\right\|_{S'(L^2;I)}$.  
	
 We start this subsection with the following remarks.
 
\begin{remark}
Since we will use the Sobolev embedding (Lemma \ref{SI}), we divide our study in three cases: $N\geq 3$ and $s<\frac{N}{2}$; $N=1,2$ and $s<\frac{N}{2}$; $N=1,2$ and $s=\frac{N}{2}$. (see respectively Lemmas \ref{LLHs1}, \ref{LLHs3} and \ref{LLHs2} bellow).     
\end{remark}  

\begin{remark}\label{Dsxb} Another interesting remark is the following claim 
\begin{equation}\label{derivadaxb}
D^s(|x|^{-b})=C_{N,b}|x|^{-b-s}.
\end{equation}
Indeed, we use the facts $\widehat{D^sf}=|\xi|^s\widehat{f}$ and $\widehat{(|x|^{-\beta})}=\frac{C_{N,\beta}}{|\xi|^{N-\beta}}$ 
for $\beta \in (0,N)$. Let $f(x)=|x|^{-b}$, we have 
$$
\widehat{D^s (|x|^{-b})}=|\xi|^s\widehat{(|x|^{-b})}=|\xi|^s\frac{C_{N,\beta}}{|\xi|^{N-b}}=\frac{C_{N,\beta}}{|\xi|^{N-(b+s)}}. 
$$
Since $0<b<\widetilde{2}$ and $0<s\leq \min\{\frac{N}{2},1\}$ then $0<b+s<N$, so taking $\beta=s+b$, we get 
$$
D^s (|x|^{-b})=\left( \frac{C_{N,\beta}}{|y|^{N-(b+s)}} \right)^\vee= C_{N,\beta}|x|^{-b-s}.
$$   
\end{remark}
\begin{lemma}\label{LLHs1} Let $N\geq 3$ and $0<b<\widetilde{2}$. If $s<\frac{N}{2}$  and $0<\alpha<\frac{4-2b}{N-2s}$ then the following statements hold  
\begin{itemize}
\item [(i)]$
\left \||x|^{-b}|u|^\alpha v \right\|_{S'(L^2;I)}\leq c (T^{\theta_1}+T^{\theta_2})\|D^s u\|^\alpha_{S(L^2;I)} \|v\|_{S(L^2;I)}
$
\item [(ii)]$
\left \|D^s(|x|^{-b}|u|^\alpha u)\right \|_{S'(L^2;I)}\leq c (T^{\theta_1}+T^{\theta_2})\|D^s u\|^{\alpha+1}_{S(L^2;I)},
$
\end{itemize}
where $I=[0,T]$ and $c,\theta_1,\theta_2 >0$.
\begin{proof} 
 (i) We divide the estimate in $B$ and $B^C$, indeed
\begin{eqnarray*}
 \hspace{1.0cm}	\left\||x|^{-b}|u|^\alpha v\right\|_{S'(L^2;I)}& \leq& \left\||x|^{-b}|u|^\alpha
v\right\|_{S'\left(L^2(B^C);I\right)}+ \left\||x|^{-b}|u|^\alpha v\right\|_{S'\left(L^2(B);I\right)}\\
&\equiv& B_1+B_2. 
\end{eqnarray*}
	
 \ First, we consider $B_1$. Let $(q_0,r_0)$ $L^2$-admissible given by\footnote{It is not difficult to check that $q_0$ and $r_0$ satisfy the conditions of admissible pair, see \eqref{L2Admissivel}.}
 \begin{equation}\label{L1Hs1} 
 q_0= \frac{4(\alpha+2)}{\alpha(N-2s)} \;\;\; \textnormal{and} \;\; r_0=\frac{ N(\alpha+2)}{ N+ \alpha s}.
 \end{equation}
 If $s<\frac{N}{2}$ then $s<\frac{N}{r_0}$ and so using the Sobolev inequality \eqref{SEI} and the H\"older inequality twice, we get 
 \begin{align}\label{L1Hs11}
 B_1\leq &\left\||x|^{-b}|u|^\alpha v\right\|_{L^{q'_0}_IL_x^{r'_0}(B^C)} \leq \left \|\||x|^{-b}\|_{L^\gamma(B^C)}\|u\|^\alpha_{L_x^{\alpha r_1}} \|v\|_{L^{r_0}_x}\right\|_{L^{q'_0}_I} \nonumber  \\
 \leq& \||x|^{-b}\|_{L^\gamma(B^C)}	\left \|  \|D^s u\|^\alpha_{L_x^{r_0}} \|v\|_{L^{r_0}_x}\right\|_{L^{q'_0}_I} \nonumber   \\
 \leq& \||x|^{-b}\|_{L^\gamma(B^C)} T^{\frac{1}{q_1}} \|D^s u\|^\alpha_{L_I^{\alpha q_2}L^{r_0}_x} \|v\|_{L^{q_0}_IL_x^{r_0}} \nonumber  \\
 =& \||x|^{-b}\|_{L^\gamma(B^C)}T^{\frac{1}{q_1}} \|D^s u\|^\alpha_{L_I^{q_0}L^{r_0}_x} \|v\|_{L^{q_0}_IL_x^{r_0}},
 \end{align}
where 
 \begin{equation}\label{CLHs1} 
\left\{\begin{array}{cl}\vspace{0.1cm}
\frac{1}{r'_0}=&\frac{1}{\gamma}+\frac{1}{r_1}+\frac{1}{r_0}\\ \vspace{0.1cm} 
\frac{1}{q'_0}=&\frac{1}{q_1}+\frac{1}{q_2}+\frac{1}{q_0}\\
q_0=& \alpha q_2\;\;,\;\;s=\frac{N}{r_0}-\frac{N}{\alpha r_1}.
\end{array}\right.
\end{equation}
In view of Remark \ref{RIxb} in order to show that the first norm in the right hand side of \eqref{L1Hs11} is bounded we need $\frac{N}{\gamma}-b<0$. Indeed, \eqref{CLHs1} is equivalent to
\begin{equation}\label{CLHs11} 
\left\{\begin{array}{cl}\vspace{0.1cm}
\frac{N}{\gamma}=&N-\frac{2N}{r_0}-\frac{N\alpha}{r_0}+\alpha s\\ \vspace{0.1cm} 
\frac{1}{q_1}=&1-\frac{\alpha+2}{q_0},
\end{array}\right.
\end{equation}
which implies, by \eqref{L1Hs1} 
\begin{equation}\label{CLHs12}
\frac{N}{\gamma}=0\;\;\;\;\;\textnormal{and}\;\;\;\;\;\frac{1}{q_1}=\frac{4-\alpha(N-2s)}{4}.
\end{equation}
 So $\frac{N}{\gamma}-b<0$ and $\frac{1}{q_1}>0$, by our hypothesis $\alpha<\frac{4-2b}{N-2s}$. Therefore, setting $\theta_1=\frac{1}{q_1}$ we deduce
\begin{equation}\label{LHsB1}
B_1 \leq c T^{\theta_1}\|D^s u\|^\alpha_{S(L^2;I)} \|v\|_{S(L^2;I)}.
\end{equation} 	

 \ We now estimate $B_2$. To do this, we use similar arguments as the ones in the estimation of $A_2$ in Lemma \ref{lemmaL2}. It follows from H\"older's inequality twice and Sobolev embedding \eqref{SEI} that
\begin{eqnarray*}
B_2&\leq &\left\||x|^{-b}|u|^\alpha v\right\|_{L^{q'}_IL_x^{r'}(B)} \leq \left \|\||x|^{-b}\|_{L^\gamma(B)} \|u\|^\alpha_{L_x^{\alpha r_1}} \|v\|_{L^r_x}\right\|_{L^{q'}_I}\\
&\leq& \left \|\||x|^{-b}\|_{L^\gamma(B)} \|D^s u\|^\alpha_{L_x^{r}} \|v\|_{L^r_x}\right\|_{L^{q'}_I}\\
&\leq&  \||x|^{-b}\|_{L^\gamma(B)}T^{\frac{1}{q_1}} \|D^s u\|^\alpha_{L_I^{\alpha q_2}L^r_x} \|v\|_{L^{q}_IL_x^r}\\
&=&  \||x|^{-b}\|_{L^\gamma(B)}T^{\frac{1}{q_1}} \|D^s u\|^\alpha_{L_I^{q}L^r_x} \|v\|_{L^{q}_IL_x^r}
\end{eqnarray*}
if $(q,r)$ $L^2$-admissible and the following system is satisfied
\begin{equation}\label{CLHs2} 
\left\{\begin{array}{cl}\vspace{0.1cm}
\frac{1}{r'}=&\frac{1}{\gamma}+\frac{1}{r_1}+\frac{1}{r}\\ \vspace{0.1cm}
s=&\frac{N}{r}-\frac{N}{\alpha r_1}, \quad s<\frac{N}{r} \\ \vspace{0.1cm}
\frac{1}{q'}=&\frac{1}{q_1}+\frac{1}{q_2}+\frac{1}{q}\\
q=& \alpha q_2 .
\end{array}\right.
\end{equation}
Similarly as in Lemma \ref{lemmaL2} we need to check that $\frac{N}{\gamma}>b$ (so that $\||x|^{-b}\|_{L^\gamma(B)}$ is finite) and $\frac{1}{q_1}>0$ for a certain choice of $(q,r)$ $L^2$-admissible pair. From \eqref{CLHs2} this is equivalent to
\begin{equation}\label{CLHs3} 
\left\{\begin{array}{cl}\vspace{0.1cm}
\frac{N}{\gamma}=& N-\frac{2N}{r}-\frac{N\alpha }{r}+\alpha s>b\\
\frac{1}{q_1}=&1-\frac{\alpha+2}{q}>0.
\end{array}\right.
\end{equation}
 The first equation in \eqref{CLHs3} implies that $\alpha<\frac{(N-b)r-2N}{N-rs}$ (assuming $s<\frac{N}{r}$), then let us choose $r$ such that 
 $$
 \frac{(N-b)r-2N}{N-rs}=\frac{4-2b}{N-2s}
 $$
 since, by our hypothesis $\alpha<\frac{4-2b}{N-2s}$. Therefore $r$ and $q$ are given by\footnote{It is easy to see that $r>2$ if, and only if, $s<\frac{N}{2}$ and $r<\frac{2N}{N-2}$ if, and only if, $b<2$. Therefore the pair $(q,r)$ given in \eqref{APLHs} is $L^2$-admissible.} 
\begin{equation}\label{APLHs}
 r=\frac{2N[N-b+2(1-s)]}{N(N-2s)+4s-bN}\;\;\textnormal{and}\;\;q=\frac{2[N-b+2(1-s)]}{N-2s},
 \end{equation}
  where we have used that $(q,r)$ is a $L^2$-admissible pair to compute the value of $q$. Note that $s<\frac{N}{r}$ if, and only if, $b+2s-N<0$. Since $s\leq 1$,  $b<\widetilde{2}$ (see \eqref{def2s}) and $N\geq 3$ it is easy to see that $s<\frac{N}{r}$ holds. In addition, from the second equation of \eqref{CLHs3} and \eqref{APLHs} we also have
\begin{equation}\label{CLHs4}
\frac{1}{q_1}=\frac{4-2b-\alpha(N-2s)}{2(N-b+2-2s)}>0,
\end{equation}
since $\alpha<\frac{4-2b}{N-2s}$.\\
Hence,
\begin{equation}\label{LHsB2} 
B_2 \leq c T^{\theta_2}\|D^s u\|^\alpha_{S(L^2;I)} \|v\|_{S(L^2;I)},
\end{equation}  
 where $\theta_2$ is given by \eqref{CLHs4}. Finally, collecting the inequalities \eqref{LHsB1} and \eqref{LHsB2} we obtain (i). 
	  
  \ (ii) Observe that
\begin{equation*}
\left\| D^s (|x|^{-b}|u|^{\alpha} u)   \right\|_{S'(L^2;I)} \leq C_1+C_2, 
 \end{equation*}
where 
\begin{align*}
 \hspace{1.2cm} C_1= &\left\| D^s (|x|^{-b}|u|^{\alpha} u)      \right\|_{S'\left(L^2(B^C);I\right)}\;\;\textnormal{and}\;\; C_2=\left\| D^s (|x|^{-b}|u|^{\alpha} u)\right\|_{S'\left(L^2(B);I\right)}.
 \end{align*}
	 
 \ To estimate $C_1$ we use the same admissible pair $(q_0,r_0)$ used to estimate the term $B_1$ in item (i). Indeed, let $$C_{11}(t)=\left\|D^s(|x|^{-b}|u|^\alpha u)\right\|_{L_x^{r'_0}(B^C)}$$ then Lemma \ref{PRODUCTRULE} (fractional product rule), Lemma \ref{CHAINRULE} (fractional chain rule) and Remark \ref{Dsxb} yield
 \begin{align}\label{L1C10}
 \hspace{0.5cm}	C_{11}(t) \leq & \||x|^{-b}\|_{L^\gamma(B^C)} \|D^s(|u|^\alpha u) \|_{L^{\beta}_x} + \|D^s(|x|^{-b})\|_{L^d(B^C)}\|u\|^{\alpha +1}_{L^{(\alpha+1)e}_x} \nonumber  \\
\leq &  \| |x|^{-b} \|_{L^\gamma(B^C)}  \|u\|^\alpha_{\alpha r_1}   \| D^s u \|_{L_x^{r_0}} + \||x|^{-b-s}\|_{L^d(B^C)}\|D^s u\|^{\alpha +1}_{L^{r_0}_x} \nonumber  \\
\leq &  \| |x|^{-b} \|_{L^\gamma(B^C)} \| D^s u \|^{\alpha+1}_{L_x^{r_0}} + \||x|^{-b-s}\|_{L^d(B^C)}\|D^s u\|^{\alpha +1}_{L^{r_0}_x},
\end{align}
where we also have used the Sobolev inequality \eqref{SEI} and \eqref{derivadaxb}. Moreover, we have the following relations
\begin{equation*}
 \left\{\begin{array}{cl}\vspace{0.1cm}
\frac{1}{r'_0}=&\frac{1}{\gamma}+\frac{1}{\beta}=\frac{1}{d}+\frac{1}{e}\\ \vspace{0.1cm}
\frac{1}{\beta}=&\frac{1}{r_1}+\frac{1}{r_0} \\ \vspace{0.1cm}
s=&\frac{N}{r_0}-\frac{N}{\alpha r_1}; \quad s<\frac{N}{r_0}\\ \vspace{0.1cm}
s=&\frac{N}{r_0}-\frac{N}{(\alpha+1)e}
\end{array}\right.
\end{equation*}
which implies that
\begin{equation}\label{L1C12} 
\left\{\begin{array}{cl}\vspace{0.1cm}
\frac{N}{\gamma}=&N-\frac{2N}{r_0}-\frac{\alpha N}{r_0}+\alpha s\\ \vspace{0.1cm}
\frac{N}{d}=&N-\frac{2N}{r_0}-\frac{\alpha N}{r_0} +\alpha s+s.
\end{array}\right.
\end{equation}
Note that, in view of \eqref{L1Hs1} we have $\frac{N}{\gamma}-b<0$ and $\frac{N}{d}-b-s<0$. These relations imply that $\| |x|^{-b} \|_{L^\gamma(B^C)}$ and $\| |x|^{-b-s} \|_{L^d(B^C)}$ are bounded quantities (see Remark \ref{RIxb}). Therefore, it follows from \eqref{L1C10} that 
\begin{equation*}
C_{11}(t)\leq  c \| D^s u \|^{\alpha+1}_{L_x^{r_0}}.
\end{equation*}
On the other hand, using $\frac{1}{q'_0}=\frac{1}{q_1}+\frac{\alpha+1}{q_0}$ and applying the H\"older inequality in the time variable we conclude
\begin{equation*}
\|C_{11}\|_{L_I^{q'_0}}\leq  c T^{\frac{1}{q_1}} \| D^s u \|^{\alpha+1}_{L_I^{q_0}L_x^{r_0}},
\end{equation*}
 where $\frac{1}{q_1}$ is given in \eqref{CLHs12}. The estimate of $C_1$ is finished since $C_1\leq\|C_{11}\|_{L_I^{q'_0}}$.
	 
\ We now consider $C_2$. Let $C_{22}(t)=\left\|D^s(|x|^{-b}|u|^\alpha u)\right\|_{L_x^{r'}(B)}$, we have $C_2\leq \|C_{22}\|_{L^{q'}_I}$. Using the same arguments as in the estimate of $C_{11}$ we obtain 
\begin{equation}\label{L1C221}
\hspace{0.5cm}C_{22}(t)\leq \| |x|^{-b} \|_{L^\gamma(B)} \| D^s u \|^{\alpha+1}_{L_x^{r}} + \||x|^{-b-s}\|_{L^d(B)}\|D^s u\|^{\alpha +1}_{L^{r}_x},
\end{equation}
 if \eqref{L1C12} is satisfied replacing $r_0$ by $r$ (to be determined later), that is
 \begin{equation}\label{L1C22} 
 \left\{\begin{array}{cl}\vspace{0.1cm}
 \frac{N}{\gamma}=&N-\frac{2N}{r}-\frac{\alpha N}{r}+\alpha s\\ \vspace{0.1cm}
 \frac{N}{d}=&N-\frac{2N}{r}-\frac{\alpha N}{r}+\alpha s+s.
\end{array}\right.
\end{equation}
In order to have that $\||x|^{-b}\|_{L^\gamma(B)}$ and $\||x|^{-b-s}\|_{L^d(B)}$ are bounded, we need $\frac{N}{\gamma}>b$ and $\frac{N}{d}>b+s$, respectively, by Remark \ref{RIxb}. Therefore, since the first equation in \eqref{L1C22} is the same as the first one in \eqref{CLHs3}, we choose $r$ as in \eqref{APLHs}. So we get $\frac{N}{\gamma}>b$, which also implies that $\frac{N}{d}-s>b$. Finally, \eqref{L1C221} and the H\"older inequality in the time variable yield
\begin{eqnarray}
C_2&\leq & cT^{\frac{1}{q_1}} \|D^s u\|^{\alpha+1}_{L_I^{(\alpha+1)q_2}L^r_x}   \nonumber    \\ 
&=& cT^{\frac{1}{q_1}} \|D^s u\|^{\alpha+1}_{L_I^{q}L^r_x},
\end{eqnarray}  
where
\begin{equation}\label{L1C222}
\frac{1}{q'}=\frac{1}{q_1}+\frac{1}{q_2}\hspace{1.0cm} q=(\alpha+1) q_2.
\end{equation}
Notice that \eqref{L1C222} is exactly to the second equation in \eqref{CLHs3}, thus $\frac{1}{q_1}>0$ (see the relation \eqref{CLHs4}). This completes the proof of Lemma \ref{LLHs1}.   
 \end{proof}
\end{lemma}
	
One important remark is that Lemma \ref{LLHs1} only holds for $N\geq 3$, since the admissible par $(q,r)$ defined in \eqref{APLHs} doesn't satisfy the condition $s<\frac{N}{r}$, for $N=1,2$. In the next lemma we study these cases.

\begin{lemma}\label{LLHs3} Let $N=1,2$ and $0<b<\widetilde{2}$. If $s<\frac{N}{2}$  and $0<\alpha<\frac{4-2b}{N-2s}$ then  
\begin{itemize}
\item [(i)]
 $
\left \||x|^{-b}|u|^\alpha v \right\|_{S'(L^2;I)}\leq c (T^{\theta_1}+T^{\theta_2})\|D^s u\|^\alpha_{S(L^2;I)} \|v\|_{S(L^2;I)}
$
\item [(ii)]
$
\left \|D^s(|x|^{-b}|u|^\alpha u)\right \|_{S'(L^2;I)}\leq c (T^{\theta_1}+T^{\theta_2})\|D^s u\|^{\alpha+1}_{S(L^2;I)},
$
\end{itemize}
 where $I=[0,T]$ and $c,\theta_1,\theta_2 >0$.
\begin{proof} 
 (i) As before, we divide the estimate in $B$ and $B^C$. The estimate on $B^C$ is the same as the term $B_1$ in Lemma \ref{LLHs1}-(i), since $(q_0,r_0)$ given in \eqref{L1Hs1} is $L^2$-admissible for $s<\frac{N}{2}$ in all dimensions. Thus we only consider the estimate on $B$.
	 
\ Indeed, set the $L^2$-admissible pair $(\bar{q},\bar{r})=(\frac{8}{2N-s},\frac{4N}{s})$. We deduce from the H\"older inequality twice and Sobolev embedding \eqref{SEI}
\begin{eqnarray*}
\left\||x|^{-b}|u|^\alpha v\right\|_{L^{\bar{q}'}_IL_x^{\bar{r}'}(B)}& \leq & \left \|\||x|^{-b}\|_{L^\gamma(B)} \|u\|^\alpha_{L_x^{\alpha r_1}} \|v\|_{L^r_x}\right\|_{L^{q'}_I}\\
&\leq&  \||x|^{-b}\|_{L^\gamma(B)}T^{\frac{1}{q_1}} \|D^s u\|^\alpha_{L_I^{\alpha q_2}L^r_x} \|v\|_{L^{q}_IL_x^r}\\
&=&  \||x|^{-b}\|_{L^\gamma(B)}T^{\frac{1}{q_1}} \|D^s u\|^\alpha_{L_I^{q}L^r_x} \|v\|_{L^{q}_IL_x^r}
\end{eqnarray*}
if $(q,r)$ is $L^2$-admissible and the following system is satisfied
\begin{equation}\label{CL3Hs2} 
\left\{\begin{array}{cl}\vspace{0.1cm}
\frac{1}{\bar{r}'}=&\frac{1}{\gamma}+\frac{1}{r_1}+\frac{1}{r}\\ \vspace{0.1cm}
s=&\frac{N}{r}-\frac{N}{\alpha r_1}; \quad s<\frac{N}{r}\\ \vspace{0.1cm}
\frac{1}{\bar{q}'}=&\frac{1}{q_1}+\frac{1}{q_2}+\frac{1}{q}\\
q=& \alpha q_2 .
\end{array}\right.
\end{equation}
Using the values of $\bar{q}$ and $\bar{r}$ given above, the previous system is equivalent to
 \begin{equation}\label{CL3Hs3} 
\left\{\begin{array}{cl}\vspace{0.1cm}
\frac{N}{\gamma}=& \frac{4(N-b)-s}{4}-\frac{N}{r}-\frac{\alpha(N-sr)}{r}+b\\
\frac{1}{q_1}=&\frac{8-2N-s}{8}-\frac{\alpha+1}{q}.
\end{array}\right.
\end{equation} 
From the first equation in \eqref{CL3Hs3} if $\alpha<\frac{r\left(4(N-b)-s\right)-4N}{N-sr}$ then $\frac{N}{\gamma}>b$, and so  $|x|^{-b}\in L^\gamma(B)$. Now, in view of the hypothesis $\alpha<\frac{4-2b}{N-2s}$ we set $r$ such that 
$$
 \frac{r\left(4(N-b)-s\right)-4N}{4(N-sr)}=\frac{4-2b}{N-2s},
$$
that is\footnote{We claim that $r$ satisfies \eqref{L2Admissivel}. In fact, obviously $r<+\infty$. Moreover $r\geq 2$ if, and only if, $8-2N+s\geq 0$ and this is true since $s>0$ and $N=1,2$.}
\begin{equation}\label{APL3Hs}
 r=\frac{4N\left(N-2s+4-2b\right)}{4s(4-2b)+(N-2s)\left(4N-4b-s\right)}.
 \end{equation}
Note that, in order to satisfy the second equation in the system \eqref{CL3Hs2} we need to verify $s<\frac{N}{r}$. A simple calculation shows that it is true if, and only if, $4b+5s<4N$ and this is true since $b<\frac{N}{3}$ and $s<\frac{N}{2}$.
	 
\ On the other hand, since we are looking for a pair $(q,r)$ $L^2$-admissible one has
\begin{equation}\label{APL3Hs1}
 q=\frac{8(N-2s+4-2b)}{(8-2N+s)(N-2s)}.
 \end{equation}
 Finally, from \eqref{APL3Hs1} the second equation in \eqref{CL3Hs3} is given by 
\begin{equation}\label{APL3Hs2}
\frac{1}{q_1}=\left( \frac{8-2N+s}{8}\right)\left(\frac{4-2b-\alpha(N-2s)}{N-2s+4-2b}\right).
\end{equation} 	  
which is positive, since $\alpha<\frac{4-2b}{N-2s}$, $s<\frac{N}{2}$ and $N=1,2$.  

 \ (ii) Similarly as in item (i) we only consider the estimate on $B$. Let $$D_2(t)=\left\||x|^{-b}|u|^\alpha u\right\|_{L_x^{\bar{r}'}(B)}.$$
  We use analogous arguments as the ones in the estimate of $C_2$ in Lemma \ref{LLHs1}-(ii). Lemmas \ref{PRODUCTRULE}-\ref{CHAINRULE}, the H\"older inequality, the Sobolev embedding \eqref{SEI} and Remark \ref{Dsxb} imply
 \begin{align}\label{L1C11}
 \hspace{0.5cm}	D_{2}(t) \leq & \||x|^{-b}\|_{L^\gamma(B)} \|D^s(|u|^\alpha u) \|_{L^{\beta}_x} + \|D^s(|x|^{-b})\|_{L^d(B)}\|u\|^{\alpha +1}_{L^{(\alpha+1)e}_x} \nonumber  \\
	   	\leq &  \| |x|^{-b} \|_{L^\gamma(B)}  \|u\|^\alpha_{\alpha r_1}   \| D^s u \|_{L_x^{r}} + \||x|^{-b-s}\|_{L^d(B)}\|D^s u\|^{\alpha +1}_{L^{r}_x} \nonumber  \\
	   	\leq &  \| |x|^{-b} \|_{L^\gamma(B)} \| D^s u \|^{\alpha+1}_{L_x^{r}} + \||x|^{-b-s}\|_{L^d(B)}\|D^s u\|^{\alpha +1}_{L^{r}_x},
	   	\end{align}
where 
\begin{equation*}
 \left\{\begin{array}{cl}\vspace{0.1cm}
\frac{1}{\bar{r}'}=&\frac{1}{\gamma}+\frac{1}{\beta}=\frac{1}{d}+\frac{1}{e}\\ \vspace{0.1cm}
\frac{1}{\beta}=&\frac{1}{r_1}+\frac{1}{r} \\ \vspace{0.1cm}
s=&\frac{N}{r}-\frac{N}{\alpha r_1}; \quad s<\frac{N}{r}\\ \vspace{0.1cm}
s=&\frac{N}{r}-\frac{N}{(\alpha+1)e},
\end{array}\right.
\end{equation*}
which is equivalent to
\begin{equation} \label{L2HsD1} 
\left\{\begin{array}{cl}\vspace{0.1cm}
 \frac{N}{\gamma}=& N-\frac{N}{\bar{r}}-\frac{(\alpha+1)N}{r}+\alpha s\\ \vspace{0.1cm}
\frac{N}{d}=&N-\frac{N}{\bar{r}}-\frac{(\alpha+1)N}{r}+\alpha s +s.
\end{array}\right.
\end{equation}
Hence, setting again $(\bar{q},\bar{r})=(\frac{8}{2N-s},\frac{4N}{s})$ the first equation in \eqref{L2HsD1} the same as the first one in \eqref{CL3Hs3}. Therefore choosing $r$ as in \eqref{APL3Hs} we have $\frac{N}{\gamma}>b$, which also implies $\frac{N}{d}>b+s$. Therefore, it follows from Remark \ref{RIxb} and \eqref{L1C11} that 
\begin{equation*}
D_2(t)\leq  c \| D^s u \|^{\alpha+1}_{L_x^{r}}.
\end{equation*}
Since, $\frac{1}{\bar{q}'}=\frac{1}{q_1}+\frac{\alpha+1}{q}$ (recall that $q$ is given in \eqref{APL3Hs1}) and applying the H\"older inequality in the time variable we conclude
\begin{equation*}
\|D_{2}\|_{L_T^{\bar{q}'}}\leq  c T^{\frac{1}{q_1}} \| D^s u \|^{\alpha+1}_{L_T^{q}L_x^{r}},
\end{equation*}
where $\frac{1}{q_1}>0$ (see \eqref{APL3Hs2}). 
\end{proof}
\end{lemma}
	
\ We finish the estimates for the nonlinearity considering the case $s=\frac{N}{2}$. Note that this case can only occur if $N=1,2$, since here we are interested in local (and global) results in $H^s(\mathbb{R}^N)$ for $\max\{0, s_c\}< s\leq \min\{\frac{N}{2},1\}$.
	
\begin{lemma}\label{LLHs2}
 Let $N=1,2$ and $0<b<\frac{N}{3}$. If $s=\frac{N}{2}$ and $0<\alpha<+\infty$ then 
\begin{itemize}
\item [(i)]$
\left\||x|^{-b}|u|^\alpha v \right\|_{S'(L^2;I)}\leq c T^{\theta_1}\| u\|^{\alpha}_{L^\infty_I H^s_x} \|v\|_{L^\infty_I L^2_x}
 $
 \item [(ii)]$
 \left\| D^s(|x|^{-b}|u|^\alpha u) \right\|_{S'(L^2;I)}\leq c  T^{\theta_1}\| u\|^{\alpha+1}_{L^\infty_I H^s_x},
$
\end{itemize}
where $I=[0,T]$ and $c,\theta_1>0$.
	 
\begin{proof}
(i) To this end we start defining the following numbers
\begin{equation}\label{L2Hs1}
 r=\frac{N(\alpha+2)}{N-2b}\;\;\; \textnormal{and} \;\;\; q=\frac{4(\alpha+2)}{N\alpha+4b},
\end{equation}
it is easy to check that $(q,r)$ is $L^2$-admissible. 

\ We divide the estimate in $B$ and $B^C$. We first consider the estimate on $B$. From H\"older's inequality 
\begin{equation}\label{LHsD21} 
\left\||x|^{-b}|u|^\alpha v\right\|_{L^{r'}_x(B)} \leq \||x|^{-b}\|_{L^\gamma(B)} \|u\|^{\alpha}_{L_x^{\alpha r_1}} \|v\|_{L^2_x},
\end{equation}
where 
\begin{equation}\label{L2Hs2}
\frac{1}{r'}=\frac{1}{\gamma}+\frac{1}{r_1}+\frac{1}{2}.
\end{equation}
In view of Remark \ref{RIxb} to show that $|x|^{-b}\in L^\gamma(B)$, we need $\frac{N}{\gamma}-b>0$. So, the relations \eqref{L2Hs1} and \eqref{L2Hs2} yield 
\begin{equation}\label{L2Hs3}
\frac{N}{\gamma}-b=\frac{\alpha(N-2b)}{2(\alpha+2)}-\frac{N}{r_1}.
\end{equation}
If we choose $\alpha r_1\in\left(\frac{2N(\alpha+2)}{N-2b},+\infty\right)$ then the right hand side of \eqref{L2Hs3} is positive. Therefore, 
\begin{equation*}
\left\||x|^{-b}|u|^\alpha v\right\|_{L^{r'}_x(B)} \leq c \|u\|^{\alpha}_{L_x^{\alpha r_1}} \|v\|_{L^2_x}.
\end{equation*}
On the other hand, since $\frac{2N(\alpha+2)}{N-2b}>2$ we can apply the Sobolev embedding  \eqref{SEI1} to obtain
\begin{equation}\label{LHsD22}
\left\||x|^{-b}|u|^\alpha v\right\|_{L^{r'}_x(B)} \leq c \|u\|^{\alpha}_{H^s} \|v\|_{L^2_x}.
\end{equation}

\ Next, we consider the estimate on $B^C$. Using the same argument as in the first case we get
\begin{equation*}
\left\||x|^{-b}|u|^\alpha v\right\|_{L^{r'}_x(B^C)} \leq \||x|^{-b}\|_{L^\gamma(B^C)} \|u\|^{\alpha}_{L_x^{\alpha r_1}} \|v\|_{L^2_x},
\end{equation*}
where the relations \eqref{L2Hs2} and \eqref{L2Hs3} hold. Thus, choosing $\alpha r_1\in\left(2,\frac{2N(\alpha+2)}{N-2b}\right)$ we have that $\frac{N}{\gamma}-b<0$, which implies $|x|^{-b}\in L^\gamma(B^C)$, by Remark \ref{RIxb}. Therefore, again the Sobolev embedding \eqref{SEI1} leads to 
\begin{equation*}
\left\||x|^{-b}|u|^\alpha v\right\|_{L^{r'}_x(B^C)} \leq c \|u\|^{\alpha}_{H_x^{s}} \|v\|_{L^2_x}.
\end{equation*}

\ Finally, it follows from the H\"older inequality in time variable, \eqref{LHsD22} and the last inequality that
\begin{equation}\label{D_2}
 \left\||x|^{-b}|u|^\alpha v\right\|_{L^{q'}_IL^{r'}_x} \leq cT^{\theta_1}\|u\|^{\alpha}_{L^\infty_I H^s} \|v\|_{L^\infty_I L^2_x},
\end{equation}
where $\theta_1=\frac{1}{q'}>0$, by \eqref{L2Hs1}. 	

 \ (ii) Similarly as in the proof of item (i) we begin setting
\begin{equation}\label{PAL2}
 r=\frac{N(\alpha+2)}{N-b-s}\;\;\;\;\textnormal{and}\;\;\;\;q=\frac{4(\alpha+2)}{\alpha N+2b+2s}.
\end{equation}
 Observe that, since $s=\frac{N}{2}$ and $0<b<\frac{N}{3}$ the denominator of $r$ is a positive number. Furthermore it is easy to verify that $(q,r)$ is $L^2$-admissible.
	 
\ First, we consider the estimate on $B$. Lemma \ref{CHAINRULE} together with the H\"older inequality and \eqref{derivadaxb} imply   
\begin{align*}
\hspace{0.5cm}	E_{1}(t) &\leq  \||x|^{-b}\|_{L^\gamma(B)} \|D^s(|u|^\alpha u) \|_{L^{\beta}_x} + \|D^s(|x|^{-b})\|_{L^d(B)}\|u\|^{\alpha +1}_{L^{(\alpha+1)e}_x} \nonumber  \\
& \leq   \| |x|^{-b} \|_{L^\gamma(B)}  \|u\|^\alpha_{L_x^{\alpha r_1}}   \| D^su \|_{L_x^{2}} + \||x|^{-b-s}\|_{L^d(B)}\| u\|^{\alpha +1}_{L^{(\alpha+1)e}_x}, 
\end{align*}
where $E_1(t)=\left\| D^s(|x|^{-b}|u|^\alpha u)\right\|_{L^{r'}_x(B)}$ and 
\begin{equation*}\label{L2E2}
\left\{\begin{array}{cl}\vspace{0.1cm}
\frac{1}{r'}=&\frac{1}{\gamma}+\frac{1}{\beta}=\frac{1}{d}+\frac{1}{e}\\ \vspace{0.1cm}
 \frac{1}{\beta}=&\frac{1}{r_1}+\frac{1}{2},
\end{array}\right.
 \end{equation*}
which implies
\begin{equation} \label{L2E3} 
\left\{\begin{array}{cl}\vspace{0.1cm}
 \frac{N}{\gamma}=& \frac{N}{2}-\frac{N}{r}-\frac{N}{r_1}\\ \vspace{0.1cm}
\frac{N}{d}=& N-\frac{N}{r}-\frac{N}{e}.
\end{array}\right.
\end{equation}

\ Now, we claim that $\| |x|^{-b} \|_{L^\gamma(B)}$ and $\| |x|^{-b-s} \|_{L^d(B)}$ are bounded quantities for a suitable choice of $r_1$ and $e$. Indeed, using the value of $r$ in  \eqref{PAL2}, \eqref{L2E3} and the fact that $s=\frac{N}{2}$ we get 
\begin{equation} \label{L2E4} 
\left\{\begin{array}{cl}\vspace{0.1cm}
\frac{N}{\gamma}-b\hspace{0.5cm}=& \frac{(\alpha+1)(N-2b)}{2(\alpha+2)}-\frac{N}{r_1}\\ \vspace{0.1cm}
\frac{N}{d}-b-s=& \frac{(\alpha+1)(N-2b)}{2(\alpha+2)}-\frac{N}{e}.
\end{array}\right.
\end{equation}
By Remark \ref{RIxb}, if $r_1, e>\frac{2N(\alpha+2)}{(\alpha+1)(N-2b)}$ then the right hand side of both equations in \eqref{L2E4} are positive, so $|x|^{-b}\in L^\gamma(B)$ and $|x|^{-b-s}\in L^d(B)$. Hence
\begin{eqnarray*}
E_{1}(t) &\leq &  c  \|u\|^\alpha_{L_x^{\alpha r_1}}   \| D^su \|_{L_x^{2}} + c\| u\|^{\alpha +1}_{L^{(\alpha+1)e}_x}.
\end{eqnarray*}
Choosing $r_1$ and $e$ as before, it is easy to see that\footnote{Increasing the value of $r_1$ if necessary.} $\alpha r_1>2$ and $(\alpha+1)e>2$, thus we can use the Sobolev inequality \eqref{SEI1}  
\begin{eqnarray}\label{L2E1}
 E_1(t) &\leq & c  \|u\|^\alpha_{H_x^{s}}   \| D^su \|_{L_x^{2}} + c\| u\|^{\alpha +1}_{H^{s}_x} \nonumber \\
 &\leq & c\|u\|^{\alpha+1}_{H_x^{s}}.
\end{eqnarray}

\ To complete the proof, we need to consider the estimate on $B^C$. By the same arguments as before we have 
\begin{align*}
\hspace{0.5cm}	E_{2}(t) & \leq   \| |x|^{-b} \|_{L^\gamma(B^C)}  \|u\|^\alpha_{L_x^{\alpha r_1}}   \| D^su \|_{L_x^{2}} + \||x|^{-b-s}\|_{L^d(B^C)}\| u\|^{\alpha +1}_{L^{(\alpha+1)e}_x}, 
\end{align*}
where $E_2(t)=\left\| D^s(|x|^{-b}|u|^\alpha u)\right\|_{L^{r'}_x(B^C)}$ and \eqref{L2E4} holds. Similarly as in item (i), since\footnote{Notice that, since $N=1,2$ and by hypothesis $\alpha>\frac{4-2b}{N}$ we have $$\frac{2N\alpha(\alpha+2)}{(\alpha+1)(N-2b)}>\frac{2N\alpha}{N-2b}>\frac{2(4-2b)}{N-2b}>2.$$} $\frac{2N\alpha(\alpha+2)}{(\alpha+1)(N-2b)}, \frac{2N(\alpha+2)}{N-b-s} >2$, we can choose $r_1$ and $e$ such that
$$
\alpha r_1\in \left(2,  \frac{2N\alpha(\alpha+2)}{(\alpha+1)(N-2b)} \right) \;\;\;\;\textnormal{and}\;\;\;\;(\alpha+1) e\in \left(2, \frac{2N(\alpha+2)}{N-2b}\right),
$$
and so we obtain from \eqref{L2E4} that $\frac{N}{\gamma}-b<0$ and $\frac{N}{d}-b-s<0$. In other words, $\| |x|^{-b} \|_{L^\gamma(B^C)}$ and $\| |x|^{-b-s} \|_{L^d(B^C)}$ are bounded quantities for these choices of $r_1$ and $e$ (see Remark \ref{RIxb}). In addition, by the Sobolev inequality \eqref{SEI1} we conclude
\begin{equation*}
 E_2(t)\leq c\|u\|^{\alpha+1}_{H_x^{s}}.
\end{equation*}

\ Finally, \eqref{L2E1} and the last inequality lead to 

 $$
\|\left\| D^s(|x|^{-b}|u|^\alpha u) \right\|_{L_I^{q'}L_x^{r'}}\leq c T^{\frac{1}{q'}}\|u\|^{\alpha+1}_{L^\infty_IH_x^{s}},
$$
where $\frac{1}{q'}>0$ by \eqref{PAL2}.
\end{proof}
\end{lemma}
	 
   \
	
  \ We now have all tools to prove the main result of this section, Theorem \ref{LWPHs}. 

\begin{proof}[\bf{Proof of Theorem \ref{LWPHs}}]	
We define $$
X= C\left([-T,T];H^s(\mathbb{R}^N)\right)\bigcap L^q\left([-T,T];H^{s,r}(\mathbb{R}^N)\right),$$ for any ($q,r$) $L^2$-admissible, and 
\begin{equation*}\label{NHs} 
\|u\|_T=\|u\|_{S\left(L^2;[-T,T]\right)}+\|D^s u\|_{S\left(L^2;[-T,T]\right)}.
\end{equation*}
We shall show that $G=G_{u_0}$ defined in \eqref{OPERATOR} is a contraction on the complete metric space
\begin{equation*}
S(a,T)=\{u \in X : \|u\|_T\leq a \}
\end{equation*}
with the metric 
$$
d_T(u,v)=\|u-v\|_{S\left(L^2;[-T,T]\right)},
$$
for a suitable choice of $a$ and $T$.

\ First, we claim that $S(a,T)$ with the metric $d_T$ is a complete metric space. Indeed, the proof follows similar arguments as in \cite{CAZENAVEBOOK} (see Theorem $1.2.5$ and the proof of Theorem 4.4.1 page 94). Since $S(a,T)\subset X$ and $X$ is a complete space, it suffices to show that $S(a,T)$, with the metric $d_T$,  is closed in $X$. Let $u_n\in S(a,T)$ such that $d_T(u_n,u)\rightarrow 0$ as $n\rightarrow +\infty$, we want to show that $u\in S(a,T)$. If $u_n\in C\left([-T,T];H^s(\mathbb{R}^N)\right)$ (see the definition of $S(a,T)$) we have, for almost all $t\in [-T,T]$, $u_n(t)$ bounded in $H^s(\mathbb{R}^N)$ and so (since $H^s(\mathbb{R}^N)$ is reflexive)
\begin{equation}\label{metricspace}
u_n(t)\rightharpoonup v(t)\;\;\textnormal{in}\;\;H^s(\mathbb{R}^N) \;\;\;\;\textnormal{and}\;\;\;\;\|v(t)\|_{H^s}\leq \liminf\limits_{n\rightarrow+\infty}\|u_n\|_{H^s}\leq a.
\end{equation}
On the other hand, the hypothesis $d_T(u_n,u)\rightarrow 0$ implies that $u_n\rightarrow u$ in $L^q_IL^r_x$ for all $(q,r)$ $L^2$-admissible. Since $(\infty,2)$ is $L^2$-admissible we get $u_n(t)\rightarrow u(t)$ in $L^2$, for almost all $t\in [-T,T]$. Therefore, by uniqueness of the limit we deduce that $u(t)=v(t)$. Also, we have from \eqref{metricspace}
$$
\|u(t)\|_{H^s}\leq a.
$$
That is, $u\in C\left([-T,T];H^s(\mathbb{R}^N)\right)$. From similar arguments, if $u_n \in L^q\left(I;H^{s,r}(\mathbb{R}^N)\right)$ we obtain $u\in S(a,I)$. This completes the proof of the claim.  

\ Returning the proof of the theorem, it follows from the Strichartz inequalities (\ref{SE1}) and \eqref{SE3} that 
		
\begin{equation}\label{NSD} 
\|G(u)\|_{S\left(L^2;[-T,T]\right)}\leq c\|u_0\|_{L^2}+ c\| F \|_{S'\left(L^2;[-T,T]\right)}
\end{equation}
and
\begin{equation}\label{NCD}  
\|D^s G(u)\|_{S\left(L^2;[-T,T]\right)}\leq c \|D^s u_0\|_{L^2}+ c\|D^s F\|_{S'\left(L^2;[-T,T]\right)},
\end{equation}
where $F(x,u)=|x|^{-b}|u|^\alpha u$. Similarly as in the proof of Theorem \ref{LWPL2}, without loss of generality we consider only the case $t > 0$. So, we deduce using Lemmas \ref{LLHs1}-\ref{LLHs3}-\ref{LLHs2} and \eqref{NHs}
\begin{equation*}\label{F} 
\| F\|_{S'(L^2;I)}\leq c(T^{\theta_1}+T^{\theta_2}) \| u \|^{\alpha+1}_I
\end{equation*} 
and
\begin{equation*}\label{DsF} 
\|D^s F\|_{S'(L^2;I)}\leq c(T^{\theta_1}+T^{\theta_2}) \| u \|^{\alpha+1}_I.
   \end{equation*}
 where $I=[0,T]$ and $\theta_1,\theta_2>0$. Hence, if $u \in S(a,T)$ then
\begin{equation*}\label{NI}
\|G(u)\|_T\leq c \|u_0\|_{H^s}+c (T^{\theta_1}+T^{\theta_2}) a^{\alpha+1}.
\end{equation*}
Now, choosing $a=2c\|u_0\|_{H^s}$ and $T>0$ such that 
\begin{equation}\label{CTHs} 
c a^{\alpha}(T^{\theta_1}+T^{\theta_2}) < \frac{1}{4},
\end{equation}
we obtain $G(u)\in S(a,T)$. Such calculations establish that $G$ is well defined on $S(a,T)$.

\ To prove that $G$ is a contraction we use \eqref{FEI} and an analogous argument as before 
\begin{eqnarray*}
d_T(G(u),G(v))&\leq &  c\| F(x,u)-F(x,v) \|_{S'\left(L^2;[-T,T]\right)}\\ 
&\leq & c (T^{\theta_1}+T^{\theta_2})\left(\|u\|^\alpha_T+\|v\|^\alpha_T\right)d_T(u,v),
\end{eqnarray*}
and so, taking $u,v\in S(a,T)$ we get
$$
 d_T(G(u),G(v))\leq c (T^{\theta_1}+T^{\theta_2})a^\alpha d_T(u,v).
$$
Therefore, from \eqref{CTHs}, $G$ is a contraction on $S(a,T)$ and by the Contraction Mapping Theorem we have a unique fixed point $u \in S(a,T)$ of $G$.
\end{proof}
		
 \ We finish this section noting that Corollary \ref{LWPH1} follows directly from Theorem \ref{LWPHs}. It is worth to mention that Corollary \ref{LWPH1} only holds for $N\geq 2$ since we assume $s\leq \min\{\frac{N}{2},1\}$ in Theorem \ref{LWPHs}.

\section{Global well-posedness}
	
\ This section is devoted to study the global well-posedness of the Cauchy problem \eqref{INLS}. Similarly as the local theory we use the fixed point theorem to prove our small data results in $H^s(\mathbb{R}^N)$. We start with a global result in $L^2(\mathbb{R}^N)$, which does not require any smallness assumption.

\subsection{$L^2$-Theory}
The global well-posedness result in $L^2(\mathbb{R}^N)$ (see Theorem \ref{GWPL2}) is an immediate consequence of Theorem \ref{LWPL2}. Indeed, using \eqref{NT3} we have that $ T(\|u_0\|_{L^2})= \frac{C}{\|u_0\|_{L^2}^d}$ for some $C,d>0$, then the conservation law \eqref{mass} allows us to reapply Theorem \ref{LWPL2} as many times as we wish preserving the length of the time interval to get a global solution. 

\subsection{$H^s$-Theory}
In this subsection, we turn our attention to proof the Theorem \ref{GWPHs}. Again the heart of the proof is to establish good estimates on the nonlinearity $F(x,u)=|x|^{-b}|u|^\alpha u$. First, we estimate the norm $\|F(x,u)\|_{S'(\dot{H}^{-s_c})}$ (see Lemma \ref{lemmaglobal1} below), next we estimate $\|F(x,u)\|_{S'(L^2)}$ (see Lemma \ref{lemmaglobal2}) and finally we consider the norm $\|D^s F(x,u)\|_{S'(L^2)}$ (see Lemmas \ref{lemmaglobal3}, \ref{lemmaglobal4} and \ref{lemmaglobal5}). 

\ We begin defining the following numbers (depending only on $N,\alpha$ and $b$)
	 
\begin{equation}\label{PHsA11}  
\widehat{q}=\frac{4\alpha(\alpha+2-\theta)}{\alpha(N\alpha+2b)-\theta(N\alpha-4+2b)}\;\;\;\widehat{r}\;=\;\frac{N\alpha(\alpha+2-\theta)}{\alpha(N-b)-\theta(2-b)}
\end{equation}
and 
 \begin{align}\label{PHsA22}
 \widetilde{a}\;=\;\frac{2\alpha(\alpha+2-\theta)}{\alpha[N(\alpha+1-\theta)-2+2b]-(4-2b)(1-\theta)}\;\;\;  \widehat{a}=\frac{2\alpha(\alpha+2-\theta)}{4-2b-(N-2)\alpha},
\end{align}
where $\theta>0$ sufficiently small\footnote{First note that, since $\theta>0$ is sufficiently small, we have that the denominators of $\widehat{q},\widehat{r},\widehat{a}$ and $\widetilde{a}$ are all positive numbers. Moreover, it is easy to see that $\widehat{r}$ satisfies \eqref{CPA2}. In fact $\widehat{a}$ can be rewritten as $\widehat{a}=\frac{\alpha+2-\theta}{1-s_c}$ and since $\theta<\alpha$ we have $\widehat{a}>\frac{2}{1-s_c}$, which implies that $\widehat{r}< \frac{2N}{N-2}$, for $N\geq 3$. We also note that $\widehat{r}\leq((\frac{2}{1-s_c})^+)'$, for $N=2$. Indeed, the last inequality is equivalent to $\varepsilon \widehat{r}<(\frac{2}{1-s_c})^+(\frac{2}{1-s_c})$ (recall \eqref{a^+}) and this is true since $\varepsilon>0$ is a small enough number. For $N=1$, we see that $\widehat{r}<\infty$. Finally, we have $\widehat{r}> \frac{2N}{N-s_c}=\frac{N\alpha}{2-b}$. Indeed, this is equivalent to $(\alpha+2-\theta)(2-b)>\alpha(N-b)-\theta(2-b) \Leftrightarrow  (\alpha+2)(2-b)>\alpha(N-b)  \Leftrightarrow  \alpha<\frac{4-2b}{N-2}$. So since $\alpha<\frac{4-2b}{N-2s}$ and $s\leq 1$ (hypothesis) we have that $\alpha<\frac{4-2b}{N-2}$ holds, consequently $\widehat{r}> \frac{2N}{N-s_c}$.}. It is easy to see that $(\widehat{q},\widehat{r})$ $L^2$-admissible, $(\widehat{a},\widehat{r})$ $\dot{H}^{s_c}$-admissible\footnote{Recall that $s_c$ is the critical Sobolev index given by $s_c=\frac{N}{2}-\frac{2-b}{\alpha}$.} and $(\widetilde{a},\widehat{r})$ $\dot{H}^{-s_c}$-admissible. Moreover, we observe that
\begin{equation}\label{EDP}  
\frac{1}{\widehat{a}}+\frac{1}{\widetilde{a}}=\frac{2}{\widehat{q}}.
\end{equation}				
		
\ Using the same notation of the previous section, we set $B=B(0,1)$ and we recall that $|x|^{-b}\in L^\gamma(B)$ if $\frac{N}{\gamma}>b$. Similarly, we have that $|x|^{-b}\in L^\gamma(B^C)$ if $\frac{N}{\gamma}<b$ (see Remark \ref{RIxb}).  
		
\ Our first result reads follows.
\begin{lemma}\label{lemmaglobal1} 
Let $ \frac{4-2b}{N}<\alpha<\alpha_s$ and $0<b<\widetilde{2}$. If $s_c<s \leq \min\{\frac{N}{2},1\}$ then the following statement holds
\begin{equation}\label{ELG11} 
\left \||x|^{-b}|u|^\alpha v \right\|_{S'(\dot{H}^{-s_c})}\leq c\| u\|^{\theta}_{L^\infty_tH^s_x}\|u\|^{\alpha-\theta}_{S(\dot{H}^{s_c})} \|v\|_{S(\dot{H}^{s_c})},
\end{equation}
where $c>0$ and $\theta\in (0,\alpha)$ is a sufficiently small number.	
	 
\begin{proof} The proof follows from similar arguments as the ones in the previous lemmas. We study the estimates in $B$ and $B^C$ separately. 

\ We first consider the set $B$. From the H\"older inequality we deduce
\begin{eqnarray}\label{LG1Hs1}
\left \|  |x|^{-b}|u|^\alpha v \right\|_{L^{\widehat{r}'}_x(B)} &\leq& \||x|^{-b}\|_{L^\gamma(B)}  \|u\|^{\theta}_{L^{\theta r_1}_x}   \|u\|^{\alpha-\theta}_{L_x^{(\alpha-\theta)r_2}}  \|v\|_{L^{\widehat{r}}_x} \nonumber  \\
&=&\||x|^{-b}\|_{L^\gamma(B)}  \|u\|^{\theta}_{L^{\theta r_1}_x}   \|u\|^{\alpha-\theta}_{L_x^{\widehat{r}}} \|v\|_{L^{\widehat{r}}_x},
\end{eqnarray}
where
\begin{equation}\label{LG1Hs2}
\frac{1}{\widehat{r}'}=\frac{1}{\gamma}+\frac{1}{r_1}+\frac{1}{r_2}+\frac{1}{\widehat{r}}\;\;\textnormal{and}\;\;\widehat{r}=(\alpha-\theta)r_2.
\end{equation}
 Now, we make use of the Sobolev embedding (Lemma \ref{SI}), so we consider two cases: $s=\frac{N}{2}$ and $s<\frac{N}{2}$.

 \ {\bf Case $s=\frac{N}{2}$}. Since $s\leq\min \{\frac{N}{2},1\}$, we only have to consider the cases where $(N,s)$ is equal to $(1,\frac{1}{2})$ or $(2,1)$. In order to have the norm $\||x|^{-b}\|_{L^\gamma(B)}$ bounded we need $\frac{N}{\gamma}>b$. In fact, observe that \eqref{LG1Hs2} implies
\begin{equation*}
\frac{N}{\gamma}=N-\frac{N(\alpha+2-\theta)}{\widehat{r}}-\frac{N}{r_1},
\end{equation*}
and from \eqref{PHsA11} it follows that
\begin{equation}\label{LG1Hs3}
\frac{N}{\gamma}-b=\frac{\theta(2-b)}{\alpha}-\frac{N}{r_1}.
\end{equation}
Since $\alpha>\frac{4-2b}{N}$ then $\frac{N\alpha}{2-b}>2$, therefore choosing
\begin{equation}\label{L1naB} 
  \theta r_1\in \left(\frac{N\alpha}{2-b},+\infty\right),
\end{equation}
 we get $\frac{N}{\gamma}>b$. Hence, inequality \eqref{LG1Hs1} and the Sobolev embedding \eqref{SEI1} yield 
\begin{equation}\label{LG1Hs4}
\left \|  |x|^{-b}|u|^\alpha v \right\|_{L^{\widehat{r}'}_x(B)} \leq c\|u\|^{\theta}_{H^s_x}  \|u\|^{\alpha-\theta}_{L_x^{\widehat{r}}} \|v\|_{L^{\widehat{r}}_x}.
\end{equation}
		
\ {\bf Case $s<\frac{N}{2}$}. Our goal here is to also obtain the inequality \eqref{LG1Hs4}. Indeed we already have the relation \eqref{LG1Hs3}, then the only change is the choice of $\theta r_1$ since we can not apply the Sobolev embedding \eqref{SEI1} when $s<\frac{N}{2}$. In this case we set 
\begin{equation}\label{LG1Hs7}
\theta r_1=\frac{2N}{N-2s},
\end{equation} 
so 
$$
\frac{N}{\gamma}-b=\theta(s-s_c)>0,
$$
that is, the quantity $\||x|^{-b}\|_{L^\gamma(B)}$ is finite. Therefore by the Sobolev embedding \eqref{SEI22} we obtain the desired inequality \eqref{LG1Hs4}. 

\ Next, we consider the set $B^C$. We claim that
\begin{equation}\label{LG1Hs41}
\left \|  |x|^{-b}|u|^\alpha v \right\|_{L^{\widehat{r}'}_x(B^C)} \leq c\|u\|^{\theta}_{H^s_x}  \|u\|^{\alpha-\theta}_{L_x^{\widehat{r}}} \|v\|_{L^{\widehat{r}}_x}.
\end{equation}
Indeed, Arguing in the same way as before we deduce
\begin{eqnarray*}
\left \|  |x|^{-b}|u|^\alpha v \right\|_{L^{\widehat{r}'}_x(B^C)} \leq \||x|^{-b}\|_{L^\gamma(B^C)}  \|u\|^{\theta}_{L^{\theta r_1}_x}   \|u\|^{\alpha-\theta}_{L_x^{\widehat{r}}} \|v\|_{L^{\widehat{r}}_x},
\end{eqnarray*}
where the relation \eqref{LG1Hs3} holds. We first show that $\||x|^{-b}\|_{L^\gamma(B^C)}$ is finite for a suitable of $r_1$. Here we also consider two cases: $s=\frac{N}{2}$ and $s<\frac{N}{2}$. In the first case, we choose $r_1$ such that
\begin{equation}\label{L1naBC} 
 \theta r_1\in\left(2,\frac{N\alpha}{2-b}\right)
\end{equation}
then, from  \eqref{LG1Hs3},  $\frac{N}{\gamma}-b<0$, so $|x|^{-b}\in L^\gamma(B^C)$. Thus, by the Sobolev inequality \eqref{SEI1} and using the last inequality we deduce \eqref{LG1Hs41}. Now if $s<\frac{N}{2}$, choosing again $\theta r_1$ as \eqref{L1naBC} one has $\frac{N}{\gamma}-b<0$. In addition, since $\alpha<\frac{4-2b}{N-2s}$ we obtain $\frac{N\alpha}{2-b}<\frac{2N}{N-2s}$, therefore the Sobolev inequality \eqref{SEI22} implies \eqref{LG1Hs41}. This completes the proof of the claim.

\ Now, inequalities \eqref{LG1Hs4} and \eqref{LG1Hs41} yield
\begin{equation}\label{LG1Hs5}
\left \|  |x|^{-b}|u|^\alpha v \right\|_{L^{\widehat{r}'}_x} \leq c\|u\|^{\theta}_{H^s_x}  \|u\|^{\alpha-\theta}_{L_x^{\widehat{r}}} \|v\|_{L^{\widehat{r}}_x}
\end{equation}
and the H\"older inequality in the time variable leads to
\begin{eqnarray*}
\left \|  |x|^{-b}|u|^\alpha v \right\|_{L_t^{\widetilde{a}'}L^{\widehat{r}'}_x}&\leq& c \|u\|^{\theta}_{L^\infty_tH^s_x} \|u\|^{\alpha-\theta}_{ L_t^{(\alpha-\theta)a_1} L_x^{\widehat{r}}} \|v\|_{L^{\widehat{a}}_tL^{\widehat{r}}_x} \nonumber  \\
&=& c \|u\|^{\theta}_{L^\infty_tH^s_x} \|u\|^{\alpha-\theta}_{ L_t^{\widehat{a}} L_x^{\widehat{r}}} \|v\|_{L^{\widehat{a}}_tL^{\widehat{r}}_x},
\end{eqnarray*}
 where
 \begin{equation*}
\frac{1}{\widetilde{a}'}=\frac{\alpha-\theta}{\widehat{a}}+\frac{1}{\widehat{a}}.
\end{equation*}
 Since $\widehat{a}$ and $\widetilde{a}$ defined in \eqref{PHsA22} satisfy the last relation we conclude the proof of \eqref{ELG11}.\footnote{Recall that $(\widehat{a},\widehat{r})$ is $\dot{H}^{s_c}$-admissible and $(\widetilde{a},\widehat{r})$ is $\dot{H}^{-s_c}$-admissible.}
\end{proof}
\end{lemma}	

 \begin{lemma}\label{lemmaglobal2} 
 Let $ \frac{4-2b}{N}<\alpha<\alpha_s$ and $0<b<\widetilde{2}$. If $s_c<s \leq \min\{\frac{N}{2},1\}$ then 
\begin{equation}\label{ELG2} 
\left\||x|^{-b}|u|^\alpha v \right\|_{S'(L^2)}\leq c\| u\|^{\theta}_{L^\infty_tH^s_x}\|u\|^{\alpha-\theta}_{S(\dot{H}^{s_c})} \| v\|_{S(L^2)},
 \end{equation}
where $c>0$ and $\theta\in (0,\alpha)$ is a sufficiently small number.
	  
\begin{proof} By the previous lemma we already have \eqref{LG1Hs5}, then applying H\"older's inequality in the time variable we obtain
 \begin{equation}\label{LGHsii}
\left\|  |x|^{-b}|u|^\alpha v\right \|_{L_t^{\widehat{q}'}L^{\widehat{r}'}_x}\leq  c \|u\|^{\theta}_{L^\infty_tH^s_x}\|u\|^{\alpha-\theta}_{L_t^{\widehat{a}} L_x^{\widehat{r}}} \|v\|_{L^{\widehat{q}}_tL^{\widehat{r}}_x},
\end{equation}
since
\begin{equation}\label{LG2Hs1}
\frac{1}{\widehat{q}'}=\frac{\alpha-\theta}{\widehat{a}}+\frac{1}{\widehat{q}}
\end{equation}
 by \eqref{PHsA11} and \eqref{PHsA22}. The proof is finished in view of $(\widehat{q},\widehat{r})$ be $L^2$-admissible.
\end{proof}
\end{lemma}
	
 \  We now estimate $\left\|D^s\left(|x|^{-b}|u|^\alpha u \right)\right\|_{S'(L^2)}$. To this end we divide our study in three cases: $N\geq 4$, $N=3$ and $N=1,2$. 
		
\begin{lemma}\label{lemmaglobal3} 
Let  $N\geq 4$, $0<b<\widetilde{2}$ and $ \frac{4-2b}{N}<\alpha<\alpha_s$. If $s_c<s \leq 1$ then the following statement holds
\begin{equation}\label{ELG3} 
\left\|D^s\left(|x|^{-b}|u|^\alpha u \right)\right\|_{S'(L^2)}\leq c\| u\|^{\theta}_{L^\infty_tH^s_x}\|u\|^{\alpha-\theta}_{S(\dot{H}^{s_c})} \| D^s u\|_{S(L^2)},
\end{equation}
where $c>0$ and $\theta\in (0,\alpha)$ is a sufficiently small number.
	      
\begin{proof} First note that we always have $s<\frac{N}{2}$ in this lemma, since we are assuming $N\geq 4$ and $s_c<s \leq 1$. Here, we also divide the estimate in  $B$ and $B^C$ separately. 

\ We begin estimating on $B$. The fractional product rule (Lemma \ref{PRODUCTRULE}) yields
\begin{eqnarray}\label{LG3Hs1} 
\left\|  D^s\left(|x|^{-b}|u|^\alpha u \right) \right \|_{L^{\widehat{r}'}_x(B)}&\leq& N_1(t,B)+N_2(t,B), 
\end{eqnarray}  
where 
$$
N_1(t,B)=\left \||x|^{-b} \right \|_{L^\gamma(B)} \left \|D^s(|u|^\alpha u)\right\|_{L^{\beta}_x}\;\;\;\;\;\;N_2(t,B)=\left\|D^s(|x|^{-b})\right\|_{L^d(B)} \left\||u|^\alpha u\right\|_{L^{e}_x}
$$    
 and
 \begin{equation}\label{LG3Hs2}
\frac{1}{\widehat{r}'}=\frac{1}{\gamma}+\frac{1}{\beta}=\frac{1}{d}+\frac{1}{e}. 
\end{equation} 
It follows from the fractional chain rule (Lemma \ref{CHAINRULE}) and H\"older's inequality that
\begin{eqnarray}\label{LG3Hs3}
 N_1(t,B) &\leq& \||x|^{-b}\|_{L^\gamma(B)}  \|u\|^{\theta}_{L^{\theta r_1}_x}   \|u\|^{\alpha-\theta}_{L_x^{(\alpha-\theta)r_2}}  \|D^s u\|_{ L^{\widehat{r}}_x}\nonumber \\
&=&   \||x|^{-b}\|_{L^\gamma(B)}   \|u\|^{\theta}_{L^{\theta r_1}_x}  \|u\|^{\alpha-\theta}_{L_x^{\widehat{r}}}  \|D^s u\|_{ L^{\widehat{r}}_x},
\end{eqnarray}
where 
\begin{equation}\label{LG3Hs4}
\frac{1}{\beta}=\frac{1}{r_1}+\frac{1}{r_2}+\frac{1}{\widehat{r}}\;\;\;\;\textnormal{and}\;\;\;\; \widehat{r}=(\alpha-\theta)r_2.
\end{equation}  
Notice that the right hand side of \eqref{LG3Hs3} is the same as the right hand side of \eqref{LG1Hs1}, with $v=D^s u$, so combining \eqref{LG3Hs2} and \eqref{LG3Hs4} we also have \eqref{LG1Hs2}. Thus, arguing in the same way as in Lemma \ref{lemmaglobal1} we obtain (recall that \eqref{LG1Hs4} also holds when $s<\frac{N}{2}$)
\begin{eqnarray}\label{LG3Hs5}
N_1(t,B) &\leq& c  \|u\|^{\theta}_{H^s_x}  \|u\|^{\alpha-\theta}_{L_x^{\widehat{r}}}  \|D^s u\|_{ L^{\widehat{r}}_x}.  
\end{eqnarray}
On the other hand, we deduce from \eqref{derivadaxb}, H\"older's inequality and the Sobolev emdebbing \eqref{SEI}
\begin{eqnarray}\label{LG3Hs7} 
N_2(t,B)&\leq& \||x|^{-b-s}\|_{L^d(B)}  \|u\|^{\theta}_{L^{\theta r_1}_x}   \|u\|^{\alpha-\theta}_{L_x^{(\alpha-\theta)r_2}}  \|u\|_{ L_x^{r_3}}\nonumber \\
&=&   \||x|^{-b-s}\|_{L^d(B)}   \|u\|^{\theta}_{L^{\theta r_1}_x} \|u\|^{\alpha-\theta}_{L^{\widehat{r}}_x}\|D^s u\|_{L^{\widehat{r}}_x},
\end{eqnarray}
where
 \begin{equation}\label{LG3Hs8}
 \left\{\begin{array}{cl}\vspace{0.1cm}
 \frac{1}{e}=&\frac{1}{r_1}+\frac{1}{r_2}+\frac{1}{r_3}\;\;\hspace{0.4cm}\;\; \widehat{r}=(\alpha-\theta)r_2\\
 s=&\frac{N}{\widehat{r}}-\frac{N}{r_3}\;\;\;\;\;\;\textnormal{with}\;\;\;\;s<\frac{N}{\widehat{r}}\;,
\end{array}\right.
\end{equation}
 which implies using \eqref{LG3Hs2} that  
\begin{equation*}
\frac{N}{d}-s= N-\frac{N(\alpha+2-\theta)}{\widehat{r}}-\frac{N}{r_1}
\end{equation*}
 and so, by \eqref{PHsA11}
\begin{equation}\label{LG3Hs9}
 \frac{N}{d}-b-s=\frac{\theta(2-b)}{\alpha}-\frac{N}{r_1}.
 \end{equation}
 Observe that the right hand side of \eqref{LG3Hs9} is the same as the right hand side of \eqref{LG1Hs3}. Hence, choosing $\theta r_1$ as in \eqref{LG1Hs7} (recall that $s<\frac{N}{2}$) we have $\frac{N}{d}-b-s>0$, so the quantity $\||x|^{-b-s}\|_{L^d(B)}$ is bounded, by Remark \ref{RIxb}. Now, the Sobolev embedding \eqref{SEI22} and \eqref{LG3Hs7} imply that
\begin{equation*}
 N_2(t,B)\leq  c \|u\|^{\theta}_{H^s_x} \|u\|^{\alpha-\theta}_{L^{\widehat{r}}_x}\|D^s u\|_{L^{\widehat{r}}_x}.
 \end{equation*}
Therefore, the last inequality together with \eqref{LG3Hs5} lead to
\begin{equation}\label{LG3Hs10}
\left\|  D^s\left(|x|^{-b}|u|^\alpha u \right) \right \|_{L^{\widehat{r}'}_x(B)}\leq c  \|u\|^{\theta}_{H^s_x} \|u\|^{\alpha-\theta}_{L^{\widehat{r}}_x}\|D^s u\|_{L^{\widehat{r}}_x}.
\end{equation}
Thus applying H\"older's inequality in the time variable and recalling \eqref{LG2Hs1},
\begin{eqnarray}\label{LG3Hs11}
\left\|  D^s\left(|x|^{-b}|u|^\alpha u \right) \right \|_{L_t^{\widehat{q}'}L^{\widehat{r}'}_x(B)} &\leq & c \|u\|^{\theta}_{L^\infty_tH^s_x} \|u\|^{\alpha-\theta}_{L^{\widehat{a}}_tL^{\widehat{r}}_x}\|D^s u\|_{L^{\widehat{q}}_tL^{\widehat{r}}_x}\nonumber \\
 &\leq&  c\| u\|^{\theta}_{L^\infty_tH^s_x}\|u\|^{\alpha-\theta}_{S(\dot{H}^{s_c})} \| D^s u\|_{S(L^2)}.
 \end{eqnarray}
 
 \ Next we consider the norm $\left\|  D^s\left(|x|^{-b}|u|^\alpha u \right) \right \|_{L^{\widehat{r}'}_x(B^C)}$. Similarly as before, replacing $B$ by $B^C$, we also get \eqref{LG3Hs3}-\eqref{LG3Hs4} and consequently by the proof of Lemma \ref{lemmaglobal1} we have the inequality \eqref{LG3Hs5}, that is
$$
N_1(t,B^C) \leq c  \|u\|^{\theta}_{H^s_x}  \|u\|^{\alpha-\theta}_{L_x^{\widehat{r}}}  \|D^s u\|_{ L^{\widehat{r}}_x}.  
$$
We also have (replacing $B$ by $B^C$)
$$
 N_2(t,B^C)\leq \||x|^{-b-s}\|_{L^d(B^C)}   \|u\|^{\theta}_{L^{\theta r_1}_x} \|u\|^{\alpha-\theta}_{L^{\widehat{r}}_x}\|D^s u\|_{L^{\widehat{r}}_x},
$$
where the relation \eqref{LG3Hs9} holds, thus setting $\theta r_1=2$ we deduce
$$
\frac{N}{d}-b-s=-\theta s_c<0,
$$
which implies that $|x|^{-b-s}\in L^d(B^C)$, by Remark \ref{RIxb}. Then the Sobolev embedding \eqref{SEI22} yields
$$
 N_2(t,B^C)\leq c   \|u\|^{\theta}_{H^s_x} \|u\|^{\alpha-\theta}_{L^{\widehat{r}}_x}\|D^s u\|_{L^{\widehat{r}}_x}.
$$
Therefore,
\begin{eqnarray*}
\left\|  D^s\left(|x|^{-b}|u|^\alpha u \right) \right \|_{L^{\widehat{r}'}_x(B^C)}&\leq &  N_1(t,B)+N_2(t,B^C)\\
&\leq & \|u\|^{\theta}_{H^s_x} \|u\|^{\alpha-\theta}_{L^{\widehat{r}}_x}\|D^s u\|_{L^{\widehat{r}}_x}.
\end{eqnarray*}

\ Finally, using H\"older's inequality in the time variable, the last inequality (recalling \eqref{LG2Hs1}) and the relation \eqref{LG3Hs11}  we get the estimate \eqref{ELG3}.  	 
 \end{proof}
\end{lemma}

 \begin{remark}
Notice that Lemma \ref{lemmaglobal3} doesn't hold in dimension three for every $\alpha < \alpha_s$ (recall \eqref{def2s}). In fact, the condition $s<\frac{N}{\widehat{r}}$ (used in \eqref{LG3Hs8}) is only true for $N\geq 4$. In the next lemma we consider the case $N=3$.
\end{remark}

\ Before stating the lemma, we define the following numbers  		
\begin{equation}\label{paradmissivel1}
 k=\frac{4\alpha(\alpha+1-\theta)}{4-2b-\alpha}\hspace{1.5cm}\;p=\frac{6\alpha(\alpha+1-\theta)}{(4-2b)(\alpha-\theta)+\alpha}
\end{equation}
 and
\begin{equation}\label{paradmissivel2}
 l=\frac{4\alpha(\alpha+1-\theta)}{\alpha(3\alpha-2+2b)-\theta(3\alpha-4+2b)},
 \end{equation}	 
where $\theta\in (0,\alpha)$. It is not difficult to verify that $(l,p)$ is $L^2$-admissible and $(k,p)$ is $\dot{H}^{s_c}$-admissible\footnote{We see that $\frac{3\alpha}{2-b}=\frac{6}{3-2s_c}<p<6$, i.e., $p$ satisfies the condition \eqref{CPA2} (and therefore \eqref{L2Admissivel}, since $\frac{6}{3-2s_c}>2$) for $N=3$.}.

\ We also define       
\begin{equation}\label{PAL41} 
 m=\frac{4D}{D-\varepsilon}\;\hspace{1.5cm}\;n=\frac{6D}{2D+\varepsilon}
\end{equation}
and 
\begin{equation}\label{PAL42}
a^*=\frac{4\theta}{2+\varepsilon-D} \;\;\;\;r^*=\frac{6\alpha\theta}{(4-2b)\theta-(2+\varepsilon-D)\alpha},
\end{equation}
where $D=\alpha-\theta+\mu$  with $\mu \in (b,1)$ and $\varepsilon$ is a sufficiently small number such that $\varepsilon<\mu-b$. Note that $2< n< 3$ ($n$ satisfies the condition \eqref{L2Admissivel} for $N=3$) and $(m,n)$ is $L^2$-admissible. Moreover, choosing $\theta=F\alpha$ with\footnote{It is easy to see that $F\in \left(\frac{1}{2},1\right)$ if $\varepsilon<\mu-b$. Therefore, since $\theta=F\alpha$, we have $\theta<\alpha$.} $F=\frac{2-\varepsilon+\mu-2b}{4-2b}$ we claim that $(a^*,r^*)$ is $\dot{H}^{s_c}$-admissible. We first show that the denominators of $a^*$ and $r^*$ are positive numbers. Indeed
$$
2+\varepsilon-D=2+\varepsilon-\mu+F\alpha-\alpha= 2+\varepsilon-\mu-\alpha(1-F)=2+\varepsilon-\mu-\alpha\left(\frac{2+\varepsilon-\mu}{4-2b}\right),
$$
so by our hypothesis $\alpha<\frac{4-2b}{3-2s}$ and since $s \leq 1$ we deduce $2+\varepsilon-D>0$. We also have (using the value of $F$ and the fact that $D>\mu$) 
$$
(4-2b)\theta-(2+\varepsilon-D)\alpha= \alpha\left((4-2b)F-2-\varepsilon+D\right)>\left(2(\mu-b)-2\varepsilon\right),
$$
which is positive setting $\varepsilon<\mu-b$.\\
Next, we show that $r^*$ satisfies the condition \eqref{CPA2}, with $N=3$. Note that $r^*$ can be rewritten as $r^*=\frac{6\alpha F}{2(\mu-b-\varepsilon)+\alpha(1-F)}$. Hence, $r^*<6$ is equivalent to 
$$
\alpha F<2(\mu-b-\varepsilon)+\alpha(1-F) \;\Leftrightarrow\;\alpha<\frac{2(\mu-b-\varepsilon)}{ 2F-1}=4-2b, 
$$
which is true since $\alpha<\frac{4-2b}{3-2s}$ and $s\leq 1$. In addition, $r^*>\frac{6}{3-2s_c}=\frac{3\alpha}{2-b}$ is equivalent to
$$
(4-2b)F>2(\mu-b-\varepsilon)+\alpha(1-F)\;\Leftrightarrow\;\alpha<4-2b.
$$
Finally, it is easy to see that $(a^*,r^*)$ satisfy the condition \eqref{CPA1}.
 
\begin{lemma}\label{lemmaglobal4} 
Let $N=3$, $ \frac{4-2b}{3}<\alpha<\frac{4-2b}{3-2s}$ and $0<b<1$. If $s_c<s \leq 1$ then there exists $\mu\in (b,1)$ such that
\begin{eqnarray}
 \left\|D^s\left(|x|^{-b}|u|^\alpha u \right)\right\|_{S'(L^2)}&\leq & c\| u\|^{\theta}_{L^\infty_tH^s_x}\|u\|^{\alpha-\theta}_{S(\dot{H}^{s_c})}\left(\| D^s u\|_{S(L^2)}+\| u\|_{S(L^2)}\right)\nonumber \\
 & & +c\| u\|^{1-\mu}_{L_t^\infty H^s_x}\|u\|^{\theta}_{S(\dot{H}^{s_c})} \|D^s u\|^{\alpha-\theta+\mu}_{S(L^2)},
\end{eqnarray}    
where $c>0$, $\theta=\alpha F$ with $F=\frac{2-\varepsilon+\mu-2b}{4-2b}$ and $\varepsilon>0$ is a sufficiently small number.
	      
\begin{proof} Observe that 
$$
\left\|D^s\left(|x|^{-b}|u|^\alpha u \right)\right\|_{S'(L^2)}\leq \left\|D^s\left(|x|^{-b}|u|^\alpha u \right)\right\|_{S'\left(L^2(B)\right)}+\left\|D^s\left(|x|^{-b}|u|^\alpha u \right)\right\|_{S'\left(L^2(B^C)\right)}.
$$
Let $A\subset \mathbb{R}^N$ that can be $B$ or $B^C$. Since $(2,6)$ is $L^2$-admissible in 3D we have
\begin{equation*}
 \left\|D^s\left(|x|^{-b}|u|^\alpha u \right)\right\|_{S'\left(L^2(A)\right)}\leq \left\|D^s\left(|x|^{-b}|u|^\alpha u \right)\right\|_{L_t^{2'}L^{6'}_x(A)}.
\end{equation*}
 As before, applying the fractional product rule (Lemma \ref{PRODUCTRULE}) we have
\begin{eqnarray}\label{LG4Hs1} 
 \left\|  D^s\left(|x|^{-b}|u|^\alpha u \right) \right \|_{L_x^{6'}(A)}&\leq& M_1(t,A)+M_2(t,A), 
\end{eqnarray}  
 where 
 $$
 M_1(t,A)=\left \||x|^{-b} \right \|_{L^\gamma(A)} \left \|D^s(|u|^\alpha u)\right\|_{L^{\beta}_x}, \quad M_2(t,A)=\left\|D^s(|x|^{-b})\right\|_{L^d(A)} \left\||u|^\alpha u\right\|_{L^{e}_x}
 $$    
and
\begin{equation}\label{LG4Hs2}
\frac{1}{6'}=\frac{1}{\gamma}+\frac{1}{\beta}=\frac{1}{d}+\frac{1}{e}. 
 \end{equation}   
	      
 \ Estimating $M_1(t,A)$. It follows by the fractional chain rule (Lemma \ref{CHAINRULE}) and H\"older's inequality that
 \begin{eqnarray}\label{LG4Hs3}
  M_1(t,A) &\leq& \||x|^{-b}\|_{L^\gamma(A)}  \|u\|^{\theta}_{L^{\theta r_1}_x}   \|u\|^{\alpha-\theta}_{L_x^{(\alpha-\theta)r_2}}  \|D^s u\|_{ L^p_x}  \nonumber \\
 &=&   \||x|^{-b}\|_{L^\gamma(A)}   \|u\|^{\theta}_{L^{\theta r_1}_x}  \|u\|^{\alpha-\theta}_{L_x^p}  \|D^s u\|_{ L^p_x},  
 \end{eqnarray}
where 
\begin{equation}\label{LG4Hs31}
 \frac{1}{\beta}=\frac{1}{r_1}+\frac{1}{r_2}+\frac{1}{p}\;\;\;\;\textnormal{and}\;\;\;\; p=(\alpha-\theta)r_2.
 \end{equation}  
  Combining \eqref{LG4Hs2} and \eqref{LG4Hs31} we obtain
 \begin{equation*}
 \frac{3}{\gamma}=\frac{5}{2}-\frac{3}{r_1}-\frac{3(\alpha+1-\theta)}{p},
\end{equation*}
 which implies, by \eqref{paradmissivel1}
\begin{equation}\label{LG4Hs32}
\frac{3}{\gamma}-b=\frac{\theta(2-b)}{\alpha}-\frac{3}{r_1}.
\end{equation}
In to order to show that $\||x|^{-b}\|_{L^\gamma(A)}$ is finite we need to verify that $\frac{3}{\gamma}-b>0$ if $A=B$ and $\frac{3}{\gamma}-b<0$ if $A=B^C$, by Remark \ref{RIxb}. Indeed if $\theta r_1=\frac{6}{3-2s}$, by \eqref{LG4Hs32}, we have 
$$\frac{3}{\gamma}-b=\theta(s-s_c)>0$$ 
and if $\theta r_1=2$ then 
$$\frac{3}{\gamma}-b=-\theta s_c<0.$$
Therefore, the inequality \eqref{LG4Hs3} and the Sobolev embedding \eqref{SEI22} yield
\begin{equation}\label{M1tA}
 M_1(t,A) \leq c \|u\|^{\theta}_{H^s_x}\|u\|^{\alpha-\theta}_{L_x^p}  \|D^s u\|_{ L^p_x}.  
\end{equation}
	     
\  We now estimate $M_2(t,A)$. Let $A=B^C$, applying the H\"older inequality and \eqref{derivadaxb} we have	
\begin{eqnarray*}
 M_2(t,B^C)& \leq & \||x|^{-b-s} \|_{L^d(B^C)} \|u\|^{\theta}_{L^{\theta r_1}_x}   \|u\|^{\alpha-\theta}_{L_x^{(\alpha-\theta)r_2}} \| u\|_{L_x^p}  \\ 
&\leq &  \||x|^{-b-s} \|_{L^d(B^C)} \|u\|^{\theta}_{L^{\theta r_1}_x}   \|u\|^{\alpha-\theta}_{L_x^p} \| u\|_{L_x^p},
\end{eqnarray*}
where 
$$
\frac{1}{e}=\frac{1}{r_1}+\frac{1}{r_2}+\frac{1}{p}\;\;\;\textnormal{and}\;\;\;p=(\alpha-\theta)r_2.
$$
 The relation \eqref{LG4Hs2} and the last relation imply
 \begin{equation*}
 \frac{3}{d}=\frac{5}{2}-\frac{3}{r_1}-\frac{3(\alpha+1-\theta)}{p}.
\end{equation*}
In view of \eqref{paradmissivel1} we deduce
\begin{equation*}
\frac{3}{d}-b=\frac{\theta(2-b)}{\alpha}-\frac{3}{r_1}.
\end{equation*}
Setting $\theta r_1=2$ we have $\frac{3}{d}-b=-\theta s_c$, so $\frac{3}{d}-b-s=-\theta s_c-s<0$, i.e., $|x|^{-b-s}\in L^d(B^C)$. So, by the Sobolev inequality \eqref{SEI22}  
\begin{eqnarray}\label{M2tBC}
 M_2(t,B^C) &\leq & c \|u\|^{\theta}_{H^s_x}   \|u\|^{\alpha-\theta}_{L_x^p} \| u\|_{L_x^p}.
\end{eqnarray}

\ We also deduce from the H\"older inequality, the Sobolev embedding\footnote{We can use the Sobolev embedding \eqref{SEI} since $s\leq 1<\frac{3}{n}$.} \eqref{SEI} and \eqref{derivadaxb}	
\begin{eqnarray*}
 M_2(t,B)& \leq& \||x|^{-b-s} \|_{L^d(B)} \|u\|^{\theta}_{L^{\theta r_1}_x}   \|u\|^{\alpha-\theta}_{L_x^{(\alpha-\theta)r_2}} \| u\|^\mu_{L_x^{\mu r_3}}\| u\|^{1-\mu}_{L_x^{(1-\mu)r_4}}  \\ 
&\leq &   \||x|^{-b-s}\|_{L^d(B)}  \|u\|^{\theta}_{L^{\theta r_1}_x} \|D^s u\|^{\alpha-\theta}_{L_x^n}  \|D^s u\|^{\mu}_{L_x^n}\| u\|^{1-\mu}_{L_x^{(1-\mu)r_4}}\\
&=& \||x|^{-b-s}\|_{L^d(B)}  \|u\|^{\theta}_{L^{r*}_x}   \|D^su\|^{\alpha-\theta+\mu}_{L_x^{n}} \|u\|^{1-\mu}_{L_x^{(1-\mu)r_4}},
\end{eqnarray*}
if the following system is satisfied
\begin{equation*}\label{LG4Hs4}
\left\{\begin{array}{cl}\vspace{0.1cm}
\frac{1}{e}=&\frac{1}{r_1}+\frac{1}{r_2}+\frac{1}{r_3}+\frac{1}{r_4}\\ \vspace{0.1cm}
 s=&\frac{3}{n}-\frac{3}{(\alpha-\theta) r_2}\;\;\;\;\;\;\;\;s=\frac{3}{n}-\frac{3}{\mu r_3}\\ \vspace{0.1cm}
r^*=&\theta r_1.
\end{array}\right.
\end{equation*}
It follows from \eqref{LG4Hs2} and the previous system that
\begin{equation}\label{LG4Hs5}
\frac{3}{d}= \frac{5}{2}+sD- \frac{3\theta}{r^*}-\frac{3D}{n}-\frac{3}{r_4},
\end{equation}
which implies by \eqref{PAL41} and \eqref{PAL42}
\begin{equation}\label{LG4Hs6}
\frac{3}{d}= \frac{7}{2}+sD- \frac{(2-b)\theta}{\alpha}-\frac{3D}{2}-\frac{3}{r_4},
\end{equation}    
where $D=\alpha-\theta+\mu$. In view of Remark \ref{RIxb} to show that $\||x|^{-b-s}\|_{L^d(B)}$ is bounded we need $\frac{3}{d}-b-s>0$. In fact, choosing $(1-\mu)r_4=\frac{6}{3-2s}$ we have 
\begin{eqnarray*}
\frac{3}{d}-b-s&=&2-b-\frac{3\alpha}{2}+\frac{3\theta}{2}+s(\alpha-\theta)-\frac{(2-b)\theta}{\alpha}\\
&=&-\alpha\left(\frac{3}{2}-\frac{2-b}{\alpha} \right)+\theta\left(\frac{3}{2}-\frac{2-b}{\alpha} \right)+s(\alpha-\theta)\\
&=&(s-s_c)(\alpha-\theta),
\end{eqnarray*}
which is positive since $s>s_c$. So $|x|^{-b-s}\in L^d(B)$ and 
\begin{equation}\label{M2tB}
 M_2(t,B)\leq c \|u\|^{1-\mu}_{H^s_x}\|u\|^{\theta}_{L^{r*}_x}   \|D^s u\|^{\alpha-\theta+\mu}_{L_x^{n}}.
\end{equation}
where we have used the Sobolev embedding \eqref{SEI22}.

\ Therefore, combining \eqref{LG4Hs1}, \eqref{M1tA} with $A=B^C$ and \eqref{M2tBC} we obtain
$$
\left\|  D^s\left(|x|^{-b}|u|^\alpha u \right) \right \|_{L_x^{6'}(B^C)}\leq c\|u\|^{\theta}_{H^s_x}   \|u\|^{\alpha-\theta}_{L_x^p} \|D^s u\|_{L_x^p}+c \|u\|^{\theta}_{H^s_x}   \|u\|^{\alpha-\theta}_{L_x^p} \| u\|_{L_x^p}.
$$
Moreover by \eqref{M1tA} with $A=B$ and \eqref{M2tB} we have
$$
\left\|  D^s\left(|x|^{-b}|u|^\alpha u \right) \right \|_{L_x^{6'}(B)}\leq c \|u\|^{\theta}_{H^s_x}   \|u\|^{\alpha-\theta}_{L_x^p} \| u\|_{L_x^p}+c \|u\|^{1-\mu}_{H^s_x}\|u\|^{\theta}_{L^{r*}_x}   \|D^s u\|^{\alpha-\theta+\mu}_{L_x^{n}}.
$$

\ Finally, since 
$$
\frac{1}{2'}=\frac{\alpha-\theta}{k}+\frac{1}{l}
$$
and 
$$
\frac{1}{2^{'}}=\frac{\theta}{a^*}+\frac{\alpha-\theta+\mu}{m},
$$
we can use H\"older's inequality in the time variable in the last two inequalities to conclude
\begin{eqnarray*}
 \left\|  D^s\left(|x|^{-b}|u|^\alpha u \right) \right \|_{L^{2'}_tL_x^{6'}(B^C)}  &\leq& c \|u\|^{\theta}_{L^\infty_tH^s_x}\|u\|^{\alpha-\theta}_{L^k_tL_x^p} \left( \|D^s u\|_{L^l_t L^p_x}+ \| u\|_{L^l_t L^p_x}\right) \\
\end{eqnarray*}
and
\begin{eqnarray*}
 \left\|  D^s\left(|x|^{-b}|u|^\alpha u \right) \right \|_{L^{2'}_tL_x^{6'}(B)}  &\leq& c \|u\|^{\theta}_{L^\infty_tH^s_x}\|u\|^{\alpha-\theta}_{L^k_tL_x^p}  \|D^s u\|_{L^l_t L^p_x}\nonumber  \\
 &&\; +\; c \|u\|^{1-\mu}_{L^\infty_tH^s_x}\|u\|^{\theta}_{L^{a^*}_tL^{r*}_x}   \|D^s u\|^{\alpha-\theta+\mu}_{L^m_tL_x^{n}},
\end{eqnarray*}

The proof is completed recalling that $(m,n)$ and $(l,p)$ are $L^2$-admissible as well as $(k,p)$ and $(a^*,r^*)$ are $\dot{H}^{s_c}$-admissible.
\end{proof}	
\end{lemma}
\begin{remark} It is worth to mention that in the previous lemma $\theta>0$ is given by $\theta=F\alpha$ and since $F<1$, we only have that $\theta<\alpha$ and it might be not true that $\theta$ is close to $0$.
\end{remark}
\ Before proving our global well-posedness result, we finish estimating the norm $\left\|D^s\left(|x|^{-b}|u|^\alpha u \right)\right\|_{S'(L^2)}$ in the dimensions $N=1,2$.  
\begin{lemma}\label{lemmaglobal5} 
Let $N=1,2$ and $\frac{4-2b}{N}<\alpha<\alpha_s$ with $0<b<\widetilde{2}$. If $s_c<s \leq \min\{\frac{N}{2},1\}$ then
\begin{eqnarray}\label{ELG5}
\left\|D^s\left(|x|^{-b}|u|^\alpha u \right)\right\|_{S'(L^2)}&\leq & c\| u\|^{\theta}_{L^\infty_tH^s_x}\|u\|^{\alpha-\theta}_{S(\dot{H}^{s_c})} \| D^s u\|_{S(L^2)}\nonumber \\
 & & +\; c\| u\|^{1+\theta}_{L_t^\infty H^s_x}\|u\|^{\alpha-\theta}_{S(\dot{H}^{s_c})},
\end{eqnarray}  
where $c>0$ and $\theta\in (0,\alpha)$ is a sufficiently small number.
\begin{proof} The proof follows from analogous arguments as the ones used in the previous lemmas. Let $A\subset \mathbb{R}^N$ that can be $B$ or $B^C$ and $(q,r)$ any $L^2$-admissible pair. By the fractional product rule (Lemma \ref{PRODUCTRULE}) we get
\begin{eqnarray}\label{LG5Hs1} 
\left\|  D^s\left(|x|^{-b}|u|^\alpha u \right) \right \|_{L^{r'}_x(A)}&\leq& P_1(t,A)+P_2(t,A), 
\end{eqnarray}  
where 
\begin{equation}\label{LG5Hs41}
 P_1(t,A)=\left \||x|^{-b} \right \|_{L^\gamma(A)} \left \|D^s(|u|^\alpha u)\right\|_{L^{\beta}_x}, \;\; P_2(t, A)=\left\|D^s(|x|^{-b})\right\|_{L^d(A)} \left\||u|^\alpha u\right\|_{L^{e}_x}
\end{equation}
 and
\begin{equation}\label{LG5Hs2}
 \frac{1}{r'}=\frac{1}{\gamma}+\frac{1}{\beta}=\frac{1}{d}+\frac{1}{e}. 
\end{equation}

 \ To estimate $P_1(t,A)$ and $P_2(t,A)$, we consider three cases: $N=1$ and $s<\frac{1}{2}$; $N=2$ and $s<1$; $N=1,2$ and $s=\frac{N}{2}$.   
	 
\ {\bf Case $N=1$ and $s<\frac{1}{2}$.}  We define the following numbers
 \begin{align}\label{paradmissivel3}  
 k^*=\frac{4\alpha(\alpha+1-\theta)}{(4-2b)(\alpha-\theta+1)-\alpha}\;\;\;\;\;\;l^*=\frac{4(\alpha+1-\theta)}{\alpha-\theta}\;\;\;\;\;\;p^*=2(\alpha+1-\theta)
\end{align}
\begin{eqnarray}\label{PAG21} 
q_0=\frac{2\alpha}{\alpha b+\theta(2-b)},\hspace{0.5cm}\textnormal{and}\hspace{0.5cm}   r_0=\frac{2\alpha}{\alpha(1-2b)-\theta(4-2b)}.
\end{eqnarray}
 It is straightforward to verify that, if $\theta>0$ is a small enough number, the assumption $0<b<\frac{1}{3}$ implies that the denominators of $q_0$, $r_0$, $k^*$ and $l^*$ are all positive numbers. Furthermore, $(q_0,r_0)$, $(l^*,p^*)$ are $L^2$-admissible\footnote{Note that, $r_0>2$ (see \eqref{L2Admissivel} for $N=1$). Moreover, since $0<b<\frac{1}{3}$ we have $p^*\geq \frac{2}{1-2s_c}=\frac{\alpha}{2-b}$ (see \eqref{CPA1} for $N=1$).} and $(k^*,p^*)$ is $\dot{H}^{s_c}$-admissible.     
	   
\ First, we estimate $P_1(t,A)$ with $r=r_0$. The fractional chain rule (Lemma \ref{CHAINRULE}) and H\"older's inequality yield
\begin{eqnarray}\label{LG5Hs3}
 P_1(t,A) &\leq& \||x|^{-b}\|_{L^\gamma(A)}  \|u\|^{\theta}_{L^{\theta r_1}_x}   \|u\|^{\alpha-\theta}_{L_x^{(\alpha-\theta)r_2}}  \|D^s u\|_{ L^{p^*}_x}  \nonumber \\
  &=&   \||x|^{-b}\|_{L^\gamma(A)}   \|u\|^{\theta}_{L^{\theta r_1}_x}  \|u\|^{\alpha-\theta}_{L_x^{p^*}}  \|D^s u\|_{ L^{p^*}_x},  
\end{eqnarray}
where 
\begin{equation}\label{LG5Hs31}
\frac{1}{\beta}=\frac{1}{r_1}+\frac{1}{r_2}+\frac{1}{p^*}\;\;\;\;\textnormal{and}\;\;\;\; p^*=(\alpha-\theta)r_2.
\end{equation}  
This implies 
\begin{equation}\label{LG5Hs32}
\frac{1}{\gamma}-b=\frac{\theta(2-b)}{\alpha}-\frac{1}{r_1},
\end{equation}  
where we have used \eqref{LG5Hs2}, \eqref{LG5Hs31}, \eqref{paradmissivel3} and \eqref{PAG21}. Now, if $A=B$ and setting $\theta r_1=\frac{2}{1-2s}$ we get $\frac{1}{\gamma}-b=\theta(s-s_c)>0$, furthermore, taking $A=B^C$ and choosing $\theta r_1=2$ one has $\frac{1}{\gamma}-b=-\theta s_c<0$. Hence, from the Sobolev embedding\footnote{Since $\theta r_1\in [2, \frac{2}{1-2s}]$ in both cases.} \eqref{SEI22} and Remark \ref{RIxb}
\begin{equation}\label{P1tA}
 P_1(t,A) \leq c \|u\|^{\theta}_{H^s_x}\|u\|^{\alpha-\theta}_{L_x^{p^*}}  \|D^s u\|_{ L^{p^*}_x}.
\end{equation}
	
\ We now consider $P_2(t,A)$ with $r=r_0$. It follows from \eqref{LG5Hs41} and \eqref{derivadaxb} that
 \begin{eqnarray}\label{LG5Hs4}
P_2(t,A)&\leq&
    \||x|^{-b-s}\|_{L^d(A)}  \|u\|^{\theta+1}_{L^{(\theta+1)e}_x}  \|u\|^{\alpha-\theta}_{L_x^\infty}
\end{eqnarray}
 and by \eqref{LG5Hs2}  
 \begin{equation}\label{LG5Hs6}
\frac{1}{d}-b=\frac{1}{2}+\frac{\theta(2-b)}{\alpha}-\frac{1}{e}.
\end{equation}
We claim that $\||x|^{-b-s}\|_{L^d(A)}$ is a finite quantity for a suitable choice of $e$. If $A=B$ we choose $(\theta+1)e=\frac{2}{1-2s}$, and if $A=B^C$ we set $(\theta+1)e=2$. We obtain in the first case  
$$
\frac{1}{d}-b-s=\theta(s-s_c)>0,
$$	   
and in the second case
 $$
 \frac{1}{\gamma}-b-s=-\theta s_c<0.
$$
So, the Sobolev embedding \eqref{SEI22}, Remark \ref{RIxb} and \eqref{LG5Hs4} yield
\begin{equation*}
P_2(t,A)\leq c\|u\|^{\theta+1}_{H^s_x}  \|u\|^{\alpha-\theta}_{L_x^\infty}.
\end{equation*}  

\ Therefore, relations \eqref{LG5Hs1}, \eqref{P1tA} and the last inequality with $A=B$ and $A=B^C$ imply that
$$
\left\|  D^s\left(|x|^{-b}|u|^\alpha u \right) \right \|_{L^{r_0'}_x(B)}\leq c \|u\|^{\theta}_{H^s_x}\|u\|^{\alpha-\theta}_{L_x^{p^*}}  \|D^s u\|_{ L^{p^*}_x}+ c\|u\|^{\theta+1}_{H^s_x}  \|u\|^{\alpha-\theta}_{L_x^\infty}
$$
and 
$$
\left\|  D^s\left(|x|^{-b}|u|^\alpha u \right) \right \|_{L^{r_0'}_x(B^C)}\leq c \|u\|^{\theta}_{H^s_x}\|u\|^{\alpha-\theta}_{L_x^{p^*}}  \|D^s u\|_{ L^{p^*}_x}+ c\|u\|^{\theta+1}_{H^s_x}  \|u\|^{\alpha-\theta}_{L_x^\infty}.
$$
Finally since
$$
\frac{1}{q_0'}=\frac{\alpha-\theta}{k^*}+\frac{1}{l^*}
$$
we apply the H\"older inequality in the time variable to get (recalling $(l^*,p^*)$ is $L^2$-admissible and $(k^*,p^*)$ is $\dot{H}^{s_c}$-admissible)
\begin{eqnarray*}
\left\|  D^s\left(|x|^{-b}|u|^\alpha u \right) \right \|_{L^{q_0'}_tL^{r_0'}_x} &\leq& c \|u\|^{\theta}_{L^\infty_tH^s_x}\|u\|^{\alpha-\theta}_{L^{k^*}_tL_x^{p^*}}  \|D^s u\|_{L^{l^*}_t L^{p^*}_x}\nonumber  \\
& &+\;c\|u\|^{\theta+1}_{L^{\infty}_t H^s_x}\|u\|^{\alpha-\theta}_{L_t^{(\alpha-\theta)q_0'}L_x^\infty}\\  
 &\leq& c\| u\|^{\theta}_{L^\infty_tH^s_x}\|u\|^{\alpha-\theta}_{S(\dot{H}^{s_c})} \| D^s u\|_{S(L^2)}\\
& &+\; c\|u\|^{\theta+1}_{L^{\infty}_t H^s_x}\|u\|^{\alpha-\theta}_{S(\dot{H}^{s_c})}.
\end{eqnarray*}
where we have used the fact that $(\alpha-\theta)q_0'=\frac{4}{1-2s_c}$, by \eqref{PAG21}, and $(\frac{4}{1-2s_c},\infty)$ is $\dot{H}^{s_c}$-admissible.
 
 \ {\bf Case $N=2$ and $s<1$.} We consider the following numbers
 \begin{align}\label{PA} 
 \widetilde{q}=\frac{2\alpha}{\alpha[b+2\varepsilon(\alpha-\theta)]+\theta(2-b)}\hspace{0.5cm}   \widetilde{r}=\frac{2\alpha}{\alpha[1-b-2\varepsilon(\alpha-\theta)]-\theta(2-b)},
 \end{align}   
\begin{eqnarray}\label{paradmissivel4}  
 l_0= \frac{ 2(\alpha+1-\theta)}{(\alpha-\theta)(1-2\varepsilon)}\;\;\hspace{0.5cm}\;\;\;p_0=\frac{2 (\alpha+1-\theta)}{1+2\varepsilon(\alpha-\theta)}
\end{eqnarray}
 and 
\begin{eqnarray}\label{paradmissivel5}  
  k_0= \frac{ 2\alpha(\alpha+1-\theta)}{\alpha[1-b-2\varepsilon(\alpha-\theta)]+(2-b)(1-\theta)}
 \end{eqnarray} 
Note that $(\widetilde{q},\widetilde{r})$, $(l_0,p_0)$ are $L^2$-admissible\footnote{The hypothesis $0<b<\frac{N}{3}$ with $N=2$ guarantee that the denominators of $\widetilde{q}$, $\widetilde{r}$, $k_0$, $l_0$ and $p_0$ are all positive numbers. Moreover, $\widetilde{r}>2$ is equivalent to $\alpha(b+2\varepsilon(\alpha-\theta))>-\theta(2-b)$ which is true, therefore $\widetilde{r}$ satisfies \eqref{L2Admissivel} for $N=2$.} and $(k_0,p_0)$ is $\dot{H}^{s_c}$-admissible\footnote{We claim that $\frac{2\alpha}{2-b}=\frac{2}{1-s_c} \leq p_0\leq((\frac{2}{1-s_c})^+)'$. Indeed, the first inequality is equivalent to $\alpha(1-b)+(1-\theta)(2-b)\geq 2\varepsilon \alpha(\alpha-\theta)$ which holds true since $\varepsilon>0$ is a small enough number. On the other hand, the later inequality holds since $\varepsilon p_0\leq (\frac{2}{1-s_c})^+(\frac{2}{1-s_c})$ (recall \eqref{a^+}) can be verified for $\varepsilon>0$ small enough.}.  

\ Estimating $P_1(t,A)$ (recall \eqref{LG5Hs41}-\eqref{LG5Hs2}) with $r=\widetilde{r}$. The fractional chain rule (Lemma \ref{CHAINRULE}) and H\"older's inequality lead to
 \begin{eqnarray}\label{LG51}
  P_1(t,A) &\leq& \||x|^{-b}\|_{L^\gamma(A)}  \|u\|^{\theta}_{L^{\theta r_1}_x}   \|u\|^{\alpha-\theta}_{L_x^{(\alpha-\theta)r_2}}  \|D^s u\|_{ L^{p_0}_x}  \nonumber \\
 &=&   \||x|^{-b}\|_{L^\gamma(A)}   \|u\|^{\theta}_{L^{\theta r_1}_x}  \|u\|^{\alpha-\theta}_{L_x^{p_0}}  \|D^s u\|_{ L^{p_0}_x},  
 \end{eqnarray}
 where
 \begin{equation}\label{LG52}
 \frac{1}{\beta}=\frac{1}{r_1}+\frac{1}{r_2}+\frac{1}{p_0}\;\;\;\;\textnormal{and}\;\;\;\; p_0=(\alpha-\theta)r_2,
 \end{equation}  
so by the relations \eqref{LG5Hs2}, \eqref{LG52}, \eqref{paradmissivel4} and \eqref{PA} one has
 \begin{equation}\label{LG53}
 \frac{2}{\gamma}-b=\frac{\theta(2-b)}{\alpha}-\frac{2}{r_1}.
 \end{equation}  
  As in the previous case, if $A=B$ we set $\theta r_1=\frac{2}{1-s}$ and then $\frac{2}{\gamma}-b>0$. On the other hand, if $A=B^C$, we set $\theta r_1=2$ and then $\frac{2}{\gamma}-b<0$. Hence, the Sobolev embedding \eqref{SEI22} and Remark \ref{RIxb} yield
\begin{equation}\label{LG54}
  P_1(t,A) \leq c \|u\|^{\theta}_{H^s_x}\|u\|^{\alpha-\theta}_{L_x^{p_0}}  \|D^s u\|_{ L^{p_0}_x}.  
\end{equation}
	   
\ Next we estimate $P_2(t,A)$ with with $r=\widetilde{r}$. An application of the H\"older inequality together with \eqref{LG5Hs41} and \eqref{derivadaxb} imply    
\begin{eqnarray*}
 P_2(t,A)  &\leq& \||x|^{-b-s}\|_{L^d(A)}\|u\|^{\theta+1}_{L^{(\theta+1)r_1}_x}  \|u\|^{\alpha-\theta}_{L_x^{(\alpha-\theta)r_2}}\\
 &\leq &\||x|^{-b-s}\|_{L^d(A)}\|u\|^{\theta+1}_{L^{(\theta+1)r_1}_x}\|u\|^{\alpha-\theta}_{L^{\frac{1}{\varepsilon}}_x},
 \end{eqnarray*}
where
\begin{equation}\label{LG5Hs8} 
 \frac{1}{e}=\frac{1}{r_1}+\frac{1}{r_2}\;,\;\;\;\;\;(\alpha-\theta)r_2=\frac{1}{\varepsilon}.
\end{equation}
 We deduce from \eqref{LG5Hs8} and \eqref{LG5Hs2} 
\begin{eqnarray}\label{LG5Hs9}
\frac{2}{d}&=& 2-\frac{2}{\widetilde{r}}-\frac{1}{r_1}-2\varepsilon(\alpha-\theta)  \nonumber  \\
 &=& 1+b+\frac{\theta(2-b)}{\alpha}-\frac{2}{r_1},
  \end{eqnarray}
where we have used \eqref{PA}. In addition, if $A=B$ and $(\theta+1)r_1=\frac{2}{1-s}$ we get
 $$
 \frac{2}{d}-b-s=\theta (s-s_c)>0,
 $$
 likewise if $A=B^C$ and $(\theta+1)r_1=2$, we have           
$$
 \frac{2}{d}-b-s=-\theta s_c-s<0.
$$
Thus
\begin{equation*}
 P_2(t,A)\leq c\|u\|^{\theta+1}_{H^s_x}\|u\|^{\alpha-\theta}_{L^{\frac{1}{\varepsilon}}_x},
\end{equation*}
where we have used the Sobolev inequality \eqref{SEI22} and and Remark \ref{RIxb}. 

\ Hence, by the relations \eqref{LG5Hs1}, \eqref{LG54} and the last inequality  
$$
\left\|  D^s\left(|x|^{-b}|u|^\alpha u \right) \right \|_{L^{\widetilde{r}}_x}\leq c \|u\|^{\theta}_{H^s_x}\|u\|^{\alpha-\theta}_{L_x^{p_0}}  \|D^s u\|_{ L^{p_0}_x}+c\|u\|^{\theta+1}_{H^s_x}\|u\|^{\alpha-\theta}_{L^{\frac{1}{\varepsilon}}_x}.
$$
Finally, from \eqref{PA} and \eqref{paradmissivel5} 
 $$
 \frac{1}{\widetilde{q}'}=\frac{\alpha-\theta}{k_0}+\frac{1}{l_0},
 $$
so applying the H\"older inequality in the time variable we deduce
 \begin{eqnarray*} 
\left\|  D^s\left(|x|^{-b}|u|^\alpha u \right) \right \|_{L^{\widetilde{q}'}_tL^{\widetilde{r}'}_x} &\leq& c \|u\|^{\theta}_{L^\infty_tH^s_x}\|u\|^{\alpha-\theta}_{L^{k_0}_tL_x^{p_0}}  \|D^s u\|_{L^{l_0}_t L^{p_0}_x} \\
& &+\;c\|u\|^{\theta+1}_{L^\infty_tH^s_x} \|u\|^{\alpha-\theta}_{{L_t^{(\alpha-\theta)\widetilde{q}'}}L^{\frac{1}{\varepsilon}}_x}\\
 &\leq& c\| u\|^{\theta}_{L^\infty_tH^s_x}\|u\|^{\alpha-\theta}_{S(\dot{H}^{s_c})} \| D^s u\|_{S(L^2)}\\
 & &+\;c\|u\|^{\theta+1}_{L^\infty_tH^s_x} \|u\|^{\alpha-\theta}_ { S(\dot{H}^{s_c})}.
\end{eqnarray*}
where we have used the fact that $(\alpha-\theta)\widetilde{q}'=\frac{2\alpha}{2-b-2\varepsilon \alpha}$ and $\left(\frac{2\alpha}{2-b-2\varepsilon \alpha},\frac{1}{\varepsilon}\right)$ is $\dot{H}^{s_c}$-admissible.
	   
 \ {\bf Case $N=1,2$ and $s=\frac{N}{2}$.} 
As before, we start defining the following numbers 
\begin{eqnarray}\label{PAG22}
\bar{a}=\frac{2(\alpha+1-\theta)}{2-s_c}\;\;\;\;\bar{q}=\frac{2(\alpha+1-\theta)}{2+s_c(\alpha-\theta)}
\end{eqnarray}
 \begin{eqnarray}\label{PAG23}
\bar{r}=\frac{2N(\alpha+1-\theta)}{N(\alpha+1-\theta)-2s_c(\alpha-\theta)-4}  
\end{eqnarray}
 and 
 \begin{align}\label{paradmissivel6}
 \bar{k}=\frac{2(\alpha+1-\theta)^2}{2(\alpha-\theta)(1-s_c)-s_c}\;\;\;\;\bar{l}=\frac{2(\alpha+1-\theta)^2}{2(\alpha-\theta)(1-s_c)+s_c\left((\alpha+1-\theta)^2-1\right)}
 \end{align}
 \begin{eqnarray}\label{paradmissivel7}
 \bar{p}=\frac{2N(\alpha+1-\theta)^2} {(N-2s_c)(\alpha+1-\theta)^2 - 4(\alpha-\theta)(1-s_c)+ 2s_c}.
\end{eqnarray}
 It is not difficult to check that $(\bar{q},\bar{r})$ and $(\bar{l},\bar{p})$ $L^2$-admissible and ($\bar{a},\bar{r}$), $(\bar{k},\bar{p})$ $\dot{H}^{s_c}$-admissible.\footnote{It is easy to see  that the denominators of $\bar{a}$ and $\bar{q}$ are positive numbers (since $s_c<1$ and $\alpha>\theta$). Furthermore, the denominators of $\bar{r}, \bar{k}, \bar{l}$ and $\bar{p}$ are also positive numbers for $\theta>0$ sufficiently small and $b<\frac{N}{3}$. We also have $\bar{r},\bar{p}\geq \frac{2N}{N-2s_c}=\frac{N\alpha}{2-b}$. Indeed $\bar{r}=\frac{2N(\alpha+1-\theta)}{N-2b-\theta(N-2s_c)}\geq \frac{N\alpha}{2-b} \Leftrightarrow \alpha(4-N)+(1-\theta)(4-2b)>-\theta\alpha(N-2s_c)$ which is true since $N=1,2$ and $\theta<1$. Moreover, $\bar{p}\geq \frac{N\alpha}{2-b}$ is equivalent to $2(\alpha-\theta)(4-2b-\alpha(N-2))\geq N\alpha-(4-2b)$ so
 $$
  \alpha\left( 2(4-2b)-N-2\alpha(N-2)  \right)+(4-2b)\geq 2\theta(4-2b-\alpha(N-2)),
 $$
 this is true since $\theta$ small enough, $N=1,2$ and $b<\frac{N}{3}$. }

 \ First, we estimate $P_1(t,A)$ with $r=\bar{r}$. The fractional chain rule (Lemma \ref{CHAINRULE}) and H\"older's inequality lead to
 \begin{eqnarray}\label{LG511}
  P_1(t,A) &\leq& \||x|^{-b}\|_{L^\gamma(A)}  \|u\|^{\theta}_{L^{\theta r_1}_x}   \|u\|^{\alpha-\theta}_{L_x^{(\alpha-\theta)r_2}}  \|D^s u\|_{ L^{\bar{p}}_x}  \nonumber \\
 &=&   \||x|^{-b}\|_{L^\gamma(A)}   \|u\|^{\theta}_{L^{\theta r_1}_x}  \|u\|^{\alpha-\theta}_{L_x^{\bar{p}}}  \|D^s u\|_{ L^{\bar{p}}_x},  
 \end{eqnarray}
 where
 \begin{equation}\label{LG522}
  \frac{1}{\beta}=\frac{1}{r_1}+\frac{1}{r_2}+\frac{1}{\bar{p}}\;\;\;\;\textnormal{and}\;\;\;\; \bar{p}=(\alpha-\theta)r_2,
 \end{equation}  
and so combining \eqref{LG5Hs2}, \eqref{LG522} \eqref{PAG23} and \eqref{paradmissivel7} we obtain
 \begin{eqnarray}\label{LG533}
\frac{N}{\gamma}-b&=&  N-b-\frac{N}{r_1}-\frac{N}{\bar{r}}-\frac{N(\alpha+1-\theta)}{\bar{p}}\nonumber\\
&=& N-b-\frac{N}{r_1}-\left( \frac{(\alpha+1-\theta)(N-2s_c)+N-2(2-s_c)}{2}  \right)\nonumber\\
&=&\frac{\theta(2-b)}{\alpha}-\frac{N}{r_1}.
\end{eqnarray}  
In order to have that the first norm in the right hand side of \eqref{LG511} is finite, we need to verify $\frac{N}{\gamma}-b>0$ if $A=B$ and $\frac{N}{\gamma}-b<0$ if $A=B^C$ for suitable choices of $r_1$. To this end, we set $r_1$ such that 
\begin{equation}\label{LG544}
\theta r_1> \frac{N\alpha}{(2-b)}\;\; (\textnormal{when } A=B) \;\;\;\textnormal{and}\;\;\; 2<\theta r_1< \frac{N\alpha}{(2-b)}\;\; (\textnormal{when } A=B^C)
\end{equation}
Hence, the Sobolev embedding \eqref{SEI1} and \eqref{LG511} yield 
\begin{eqnarray}\label{LG555}
 P_1(t,A) &\leq& c \|u\|^{\theta}_{H^s_x}  \|u\|^{\alpha-\theta}_{L_x^{\bar{p}}}  \|D^s u\|_{ L^{\bar{p}}_x}.
\end{eqnarray}

\ We now consider $P_2(t,A)$ with $r=\bar{r}$. By the H\"older inequality and \eqref{LG5Hs41} 
 \begin{eqnarray}\label{LG5Hs11}
   P_2(t,A) &\leq& \||x|^{-b-s}\|_{L^d(A)}\|u\|^{\theta+1}_{L^{(\theta+1)r_1}_x}  \|u\|^{\alpha-\theta}_{L_x^{(\alpha-\theta)r_2}} \nonumber  \\
  &=&\||x|^{-b-s}\|_{L^d(A)}\|u\|^{\theta+1}_{L^{(\theta+1)r_1}_x}\|u\|^{\alpha-\theta}_{L_x^{\bar{r}}},
\end{eqnarray}
 where
 \begin{equation}\label{LG5Hs12}
  \frac{1}{e}=\frac{1}{r_1}+\frac{1}{r_2}\;\;\;\;\textnormal{and}\;\;\;\;\bar{r}=(\alpha-\theta)r_2.
 \end{equation}
The relations \eqref{LG5Hs2} and \eqref{LG5Hs12} as well as $\bar{r}$ defined in \eqref{PAG23}, yield (recall $s=\frac{N}{2}$)  
\begin{eqnarray}\label{LG5Hs13}
\frac{N}{d}-b-s &=& N-b-s-\frac{N}{r_1}-\frac{N(\alpha+1-\theta)}{\bar{r}} \nonumber \\
&=&  \frac{N}{2}+(2-b)-\frac{N}{r_1}-\frac{N(\alpha+1-\theta)}{2}+s_c(\alpha-\theta) \nonumber\\
&=&\frac{\theta(2-b)}{\alpha}-\frac{N}{r_1}.
\end{eqnarray}   	   
We see that the right hand side of \eqref{LG5Hs13} is equal to the right hand side of \eqref{LG533}, so choosing $r_1$ as in \eqref{LG544} and again applying the Sobolev inequality \eqref{SEI1}, we conclude  
\begin{equation*}
 P_2(t,A)\leq c \|u\|^{\theta+1}_{H^s_x}\|u\|^{\alpha-\theta}_{L_x^{\bar{r}}}.
\end{equation*}

\ The inequalities \eqref{LG5Hs1}, \eqref{LG555} and the last inequality imply that
$$
\left\|  D^s\left(|x|^{-b}|u|^\alpha u \right) \right \|_{L^{\bar{r}'}_x}\leq c \|u\|^{\theta}_{H^s_x}  \|u\|^{\alpha-\theta}_{L_x^{\bar{p}}}  \|D^s u\|_{ L^{\bar{p}}_x}+c \|u\|^{\theta+1}_{H^s_x}\|u\|^{\alpha-\theta}_{L_x^{\bar{r}}}.
$$
Since
\begin{equation*}
\frac{1}{\bar{q}'}=\frac{\alpha-\theta}{\bar{k}}+\frac{1}{\bar{l}}
\end{equation*}
we can apply the H\"older inequality in the time variable to deduce
\begin{eqnarray*}
\left\|  D^s\left(|x|^{-b}|u|^\alpha u \right) \right \|_{L^{\bar{q}'}_tL^{\bar{r}'}_x} &\leq& c \|u\|^{\theta}_{L^\infty_tH^s_x}\|u\|^{\alpha-\theta}_{L^{\bar{k}}_tL_x^{\bar{p}}}  \|D^s u\|_{L^{\bar{l}}_t L^{\bar{p}}_x}\\
&&+c \;\|u\|^{\theta+1}_{L^{\infty}_t H^s_x}\|u\|^{\alpha-\theta}_{L_t^{(\alpha-\theta)\bar{q}'}L_x^{\bar{r}}}\\
&\leq& c\| u\|^{\theta}_{L^\infty_tH^s_x}\|u\|^{\alpha-\theta}_{S(\dot{H}^{s_c})} \| D^s u\|_{S(L^2)}\\
&& +c\|u\|^{\theta+1}_{L^{\infty}_t H^s_x}\|u\|^{\alpha-\theta}_{L_t^{\bar{a}}L_x^{\bar{r}}},
\end{eqnarray*}
where in the last equality we have used the fact that $\bar{a}=(\alpha-\theta)\bar{q}'$. This  completes the proof since $(\bar{a},\bar{r})$ $\dot{H}^{s_c}$-admissible. 
 \end{proof}	
 \end{lemma}

\ The next result follows directly from Lemmas \ref{lemmaglobal3}, \ref{lemmaglobal4} and \ref{lemmaglobal5}.

\begin{corollary}\label{corollarylemmas} Assume $ \frac{4-2b}{N}<\alpha<\alpha_s$ and $0<b<\widetilde{2}$. If $s_c<s\leq \min\{\frac{N}{2}, 1\}$ then following statement hold: 
\begin{eqnarray*}\label{corollaryglobal}
\|D^s F\|_{S'(L^2)}\;\;&\leq& c\| u\|^{\theta}_{L^\infty_tH^s_x}\|u\|^{\alpha-\theta}_{S(\dot{H}^{s_c})}\left( \|D^s u\|_{S(L^2)}+\| u\|_{S(L^2)}+\| u\|_{L^\infty_tH^s_x}\right)\\
& &+\; c\| u\|^{1-\mu}_{L^\infty_tH^s_x}\|u\|^{\theta}_{S(\dot{H}^{s_c})} \|D^s u\|^{\alpha-\theta+\mu}_{S(L^2)},
\end{eqnarray*}
where $F(x,u)=|x|^{-b}|u|^\alpha u$.
\end{corollary}
	
\ Now, we have all the tools to prove the Theorem \ref{GWPHs}. Similarly as in the local theory, we use the contraction mapping principle. 
		
\begin{proof}[\bf{Proof of Theorem \ref{GWPHs}}] First, we define
$$
B=\{ u:\;\|u\|_{S(\dot{H}^{s_c})}\leq 2\|U(t)u_0\|_{S(\dot{H}^{s_c})}\;\textnormal{and}\;\|u\|_{S(L^2)}+\|D^s u\|_{S(L^2)}\leq 2c\|u_0\|_{H^s}\}.
$$ We prove that $G=G_{u_0}$ defined in \eqref{OPERATOR} is a contraction on $B$ equipped with the metric 
$$
d(u,v)=\|u-v\|_{S(L^2)}+\|u-v\|_{S(\dot{H}^{s_c})}.
$$
		
\ Indeed, we deduce by the Strichartz inequalities (\ref{SE1}), (\ref{SE2}), \eqref{SE3} and \eqref{SE5}
\begin{equation}\label{GHs1}
\|G(u)\|_{S(\dot{H}^{s_c})}\leq \|U(t)u_0\|_{S(\dot{H}^{s_c})}+ c\| F \|_{S'(\dot{H}^{-s_c})}
\end{equation}
\begin{equation}\label{GHs2}
\|G(u)\|_{S(L^2)}\leq c\|u_0\|_{L^2}+ c\| F \|_{S'(L^2)}
\end{equation}
and	
\begin{equation}\label{GHs3}
 \|D^s G(u)\|_{S(L^2)}\leq c \|D^s u_0\|_{L^2}+ c\|D^s F\|_{S'(L^2)},
\end{equation}
where $F(x,u)=|x|^{-b}|u|^\alpha u$. On the other hand, it follows from Lemmas \ref{lemmaglobal1} and \ref{lemmaglobal2} together with Corollary \ref{corollarylemmas} that
\begin{eqnarray*}
\|F\|_{S'(\dot{H}^{-s_c})}&\leq & c\| u \|^\theta_{L^\infty_tH^s_x}\| u \|^{\alpha-\theta}_{S(\dot{H}^{s_c})}\| u \|_{S(\dot{H}^{s_c})}\\
\|F\|_{S'(L^2)}&\leq& c\| u \|^\theta_{L^\infty_tH^s_x}\| u \|^{\alpha-\theta}_{S(\dot{H}^{s_c})}\| u \|_{S(L^2)}
\end{eqnarray*}
and	
\begin{eqnarray*}
\|D^s F\|_{S'(L^2)}&\leq& c\| u\|^{\theta}_{L^\infty_tH^s_x}\|u\|^{\alpha-\theta}_{S(\dot{H}^{s_c})}\left( \|D^s u\|_{S(L^2)}+\| u\|_{S(L^2)}+\| u\|_{L^\infty_tH^s_x}\right)\\
& & +\; c\| u\|^{1-\mu}_{L^\infty_tH^s_x}\|u\|^{\theta}_{S(\dot{H}^{s_c})} \|D^s u\|^{\alpha-\theta+\mu}_{S(L^2)}.
\end{eqnarray*}
Combining \eqref{GHs1}-\eqref{GHs3} and the last inequalities, we get for $u\in B$
\begin{align*}\label{TGHS}
\|G(u)\|_{S(\dot{H}^{s_c})}\leq& \|U(t)u_0\|_{S(\dot{H}^{s_c})} +c\| u \|^\theta_{L^\infty_tH^s_x}\| u \|^{\alpha-\theta}_{S(\dot{H}^{s_c})}\| u \|_{S(\dot{H}^{s_c})} \nonumber \\
\leq & \|U(t)u_0\|_{S(\dot{H}^{s_c})}+2^{\alpha+1}c^{\theta+1}\|u_0\|^\theta_{H^s}\| U(t)u_0 \|^{\alpha-\theta+1}_{S(\dot{H}^{s_c})}.
\end{align*}
In addition, setting $X=\|D^s u\|_{S(L^2)}+\| u\|_{S(L^2)}+\| u\|_{L^\infty_tH^s_x}$ 
\begin{eqnarray*}
\|G(u)\|_{S(L^2)}+\|D^s G(u)\|_{S(L^2)}&\leq & c\|u_0\|_{H^s}+c\| u \|^\theta_{L^\infty_tH^s_x}\| u \|^{\alpha-\theta}_{S(\dot{H}^{s_c})}X \nonumber\\
& & +\;c\| u\|^{1-\mu}_{L^\infty_tH^s_x}\|u\|^{\theta}_{S(\dot{H}^{s_c})} \|D^s u\|^{\alpha-\theta+\mu}_{S(L^2)}\\
&\leq & c\|u_0\|_{H^s}+2^{\alpha+2}c^{\theta+2}\|u_0\|_{H^s}^{\theta+1}\| U(t)u_0 \|^{\alpha-\theta}_{S(\dot{H}^{s_c})}\\
& & + 2^{\alpha+1}c^{\alpha-\theta+2}\|u_0\|_{H^s}^{\alpha-\theta+1}\| U(t)u_0 \|^{\theta}_{S(\dot{H}^{s_c})},
\end{eqnarray*}
where we have have used the fact that $X\leq 2^2c\|u_0\|_{H^s}$ since $u\in B$.

\ Now if $\| U(t)u_0 \|_{S(\dot{H}^{s_c})}<\delta$ with 
\begin{equation}\label{WD}
\delta\leq \min\left\{\sqrt[\alpha-\theta]{\frac{1}{2c^{\theta+1}2^{\alpha+1}A^\theta}}     , \sqrt[\alpha-\theta]{ \frac{1}{4c^{\theta+1}2^{\alpha+2}A^\theta}},\sqrt[\theta]{ \frac{1}{4c^{\alpha-\theta+1}2^{\alpha+1}A^{\alpha-\theta}}}\right\},
\end{equation}
where $A>0$ is a number such that $\|u_0\|_{H^s}\leq A$, we get 
$$\|G(u)\|_{S(\dot{H}^{s_c})}\leq 2\|U(t)u_0\|_{S(\dot{H}^{s_c})}$$ and  $$\|G(u)\|_{S(L^2)}+\|D^s G(u)\|_{S(L^2)}\leq 2c\|u_0\|_{H^s},$$ that is  $G(u)\in B$.

\ To complete the proof we show that $G$ is a contraction on $B$. From \eqref{FEI} and repeating the above computations one has
\begin{align*}\label{C1GH1}
\|G(u)-G(v)\|_{S(\dot{H}^{s_c})}\leq&c\|F(x,u)-F(x,v)\|_{S(\dot{H}^{-s_c})}\\
\leq& c\left\||x|^{-b}|u|^{\alpha}|u-v|\right\|_{S(\dot{H}^{-s_c})}+\left\||x|^{-b}|v|^{\alpha}|u-v|\right\|_{S(\dot{H}^{-s_c})}\\
\leq & c\| u \|^\theta_{L^\infty_tH^s_x}\| u \|^{\alpha-\theta}_{S(\dot{H}^{s_c})}\| u -v\|_{S(\dot{H}^{s_c})}\\
 &+c \| v \|^\theta_{L^\infty_tH^s_x}\| v \|^{\alpha-\theta}_{S(\dot{H}^{s_c})}\| u -v\|_{S(\dot{H}^{s_c})}
\end{align*}
which implies, taking $u,v\in B$
\begin{eqnarray*}\label{C1GH1}
\|G(u)-G(v)\|_{S(\dot{H}^{s_c})} &\leq&  2c(2c)^\theta \| u_0 \|^\theta_{H^s}2^{\alpha-\theta}\|U(t)u_0 \|^{\alpha-\theta}_{S(\dot{H}^{s_c})}\| u -v\|_{S(\dot{H}^{s_c})}\\
& =&   2^{\alpha+1}c^{\theta+1} \| u_0 \|^\theta_{H^s}\|U(t)u_0 \|^{\alpha-\theta}_{S(\dot{H}^{s_c})}\| u -v\|_{S(\dot{H}^{s_c})}.
\end{eqnarray*}
By similar arguments we also obtain
\begin{equation*}\label{C2GH1}
\|G(u)-G(v)\|_{S(L^2)}\leq 2^{\alpha+1}c^{\theta+1} \| u_0 \|^\theta_{H^s}\|U(t)u_0 \|^{\alpha-\theta}_{S(\dot{H}^{s_c})} \|u-v\|_{S(L^2)}.
\end{equation*}
Finally, from the two last inequalities and \eqref{WD} 
$$
d(G(u),G(v)) \leq 2^{\alpha+1}c^{\theta+1} \| u_0 \|^\theta_{H^s}\|U(t)u_0 \|^{\alpha-\theta}_{S(\dot{H}^{s_c})}d(u,v)\leq \frac{1}{2}d(u,v),
$$
i.e., $G$ is a contraction.

\ Therefore, by the Banach Fixed Point Theorem, $G$ has a unique fixed point $u\in B$, which is a global solution of \eqref{INLS}.
\end{proof}
\subsection*{Acknowledgments} 

The author would like to thank professor L. Farah (UFMG) for the very useful suggestions given during the elaboration of this paper. The author was supported by CAPES/Brazil.
 
\bibliographystyle{abbrv}
\bibliography{bibguzman}	

\begin{thebibliography}{10}

\bibitem{BERLOF}
J.~Bergh and J.~L{\"o}fstr{\"o}m.
\newblock {\em Interpolation spaces. {A}n introduction}.
\newblock Springer-Verlag, Berlin-New York, 1976.
\newblock Grundlehren der Mathematischen Wissenschaften, No. 223.

\bibitem{CAZENAVEBOOK}
T.~Cazenave.
\newblock {\em Semilinear {S}chr\"odinger equations}, volume~10 of {\em Courant
  Lecture Notes in Mathematics}.
\newblock New York University, Courant Institute of Mathematical Sciences, New
  York; American Mathematical Society, Providence, RI, 2003.

\bibitem{CAZENAVECONTINUOUS}
T.~Cazenave, D.~Fang, and Z.~Han.
\newblock Continuous dependence for {NLS} in fractional order spaces.
\newblock {\em Ann. Inst. H. Poincar\'e Anal. Non Lin\'eaire}, 28(1):135--147,
  2011.

\bibitem{CAZEWE}
T.~Cazenave and F.~B. Weissler.
\newblock Some remarks on the nonlinear {S}chr\"odinger equation in the
  critical case.
\newblock In {\em Nonlinear semigroups, partial differential equations and
  attractors}, volume 1394 of {\em Lecture Notes in Math.}, pages 18--29.
  Springer, Berlin, 1989.

\bibitem{CAZENAVEHS}
T.~Cazenave and F.~B. Weissler.
\newblock The {C}auchy problem for the critical nonlinear {S}chr\"odinger
  equation in {$H\sp s$}.
\newblock {\em Nonlinear Anal.}, 14(10):807--836, 1990.

\bibitem{CAZWEI}
T.~Cazenave and F.~B. Weissler.
\newblock Rapidly decaying solutions of the nonlinear {S}chr\"odinger equation.
\newblock {\em Comm. Math. Phys.}, 147(1):75--100, 1992.

\bibitem{CHRWEI}
F.~M. Christ and M.~I. Weinstein.
\newblock Dispersion of small amplitude solutions of the generalized
  {K}orteweg-de {V}ries equation.
\newblock {\em J. Funct. Anal.}, 100(1):87--109, 1991.

\bibitem{DEMENGEL}
F.~Demengel and G.~Demengel.
\newblock {\em Functional spaces for the theory of elliptic partial
  differential equations}.
\newblock Universitext. Springer, London; EDP Sciences, Les Ulis, 2012.

\bibitem{DUCHOLROU}
T.~Duyckaerts, J.~Holmer, and S.~Roudenko.
\newblock Scattering for the non-radial 3{D} cubic nonlinear {S}chr\"odinger
  equation.
\newblock {\em Math. Res. Lett.}, 15(6):1233--1250, 2008.

\bibitem{LG}
L.~G. Farah.
\newblock Global well-posedness and blow-up on the energy space for the
  inhomogeneous nonlinear {S}chr\"odinger equation.
\newblock {\em J. Evol. Equ.}, 16(1):193--208, 2016.

\bibitem{FIBI}
G.~Fibich and X.~P. Wang.
\newblock Stability of solitary waves for nonlinear {S}chr\"odinger equations
  with inhomogeneous nonlinearities.
\newblock {\em Phys. D}, 175(1-2):96--108, 2003.

\bibitem{GENOUD}
F.~Genoud.
\newblock An inhomogeneous, {$L\sp 2$}-critical, nonlinear {S}chr\"odinger
  equation.
\newblock {\em Z. Anal. Anwend.}, 31(3):283--290, 2012.

\bibitem{GENSTU}
F.~Genoud and C.~A. Stuart.
\newblock Schr\"odinger equations with a spatially decaying nonlinearity:
  existence and stability of standing waves.
\newblock {\em Discrete Contin. Dyn. Syst.}, 21(1):137--186, 2008.

\bibitem{GILL}
T.~S. Gill.
\newblock Optical guiding of laser beam in nonuniform plasma.
\newblock {\em Pramana J. Phys}, 55(5-6):835--842, 2000.

\bibitem{GV}
J.~Ginibre and G.~Velo.
\newblock On a class of nonlinear {S}chr\"odinger equations. {I}. {T}he
  {C}auchy problem, general case.
\newblock {\em J. Funct. Anal.}, 32(1):1--32, 1979.

\bibitem{GV1}
J.~Ginibre and G.~Velo.
\newblock The global {C}auchy problem for the nonlinear {S}chr\"odinger
  equation revisited.
\newblock {\em Ann. Inst. H. Poincar\'e Anal. Non Lin\'eaire}, 2(4):309--327,
  1985.

\bibitem{GUEVARA}
C.~D. Guevara.
\newblock Global behavior of finite energy solutions to the {$d$}-dimensional
  focusing nonlinear {S}chr\"odinger equation.
\newblock {\em Appl. Math. Res. Express. AMRX}, 2014(2):177--243, 2014.

\bibitem{HOLROU}
J.~Holmer and S.~Roudenko.
\newblock A sharp condition for scattering of the radial 3{D} cubic nonlinear
  {S}chr\"odinger equation.
\newblock {\em Comm. Math. Phys.}, 282(2):435--467, 2008.

\bibitem{KATO}
T.~Kato.
\newblock On nonlinear {S}chr\"odinger equations.
\newblock {\em Ann. Inst. H. Poincar\'e Phys. Th\'eor.}, 46(1):113--129, 1987.

\bibitem{tao:keel}
M.~Keel and T.~Tao.
\newblock Endpoint {S}trichartz estimates.
\newblock {\em Amer. J. Math.}, 120(5):955--980, 1998.

\bibitem{KENIG}
C.~E. Kenig and F.~Merle.
\newblock Global well-posedness, scattering and blow-up for the
  energy-critical, focusing, non-linear {S}chr\"odinger equation in the radial
  case.
\newblock {\em Invent. Math.}, 166(3):645--675, 2006.

\bibitem{FELGUS}
F.~Linares and G.~Ponce.
\newblock {\em Introduction to nonlinear dispersive equations}.
\newblock Universitext. Springer, New York, second edition, 2015.

\bibitem{LIU}
C.~S. Liu and V.~K. Tripathi.
\newblock Laser guiding in an axially nonuniform plasma channel.
\newblock {\em Physics of Plasmas}, 1(9):3100--3103, 1994.

\bibitem{MERLE}
F.~Merle.
\newblock Nonexistence of minimal blow-up solutions of equations {$iu\sb
  t=-\Delta u-k(x)\vert u\vert \sp {4/N}u$} in {${\bf R}\sp N$}.
\newblock {\em Ann. Inst. H. Poincar\'e Phys. Th\'eor.}, 64(1):33--85, 1996.

\bibitem{RAFAEL}
P.~Rapha{\"e}l and J.~Szeftel.
\newblock Existence and uniqueness of minimal blow-up solutions to an
  inhomogeneous mass critical {NLS}.
\newblock {\em J. Amer. Math. Soc.}, 24(2):471--546, 2011.

\bibitem{TSUTSUMI}
Y.~Tsutsumi.
\newblock {$L\sp 2$}-solutions for nonlinear {S}chr\"odinger equations and
  nonlinear groups.
\newblock {\em Funkcial. Ekvac.}, 30(1):115--125, 1987.

\bibitem{WEINSTEIN}
M.~I. Weinstein.
\newblock Nonlinear {S}chr\"odinger equations and sharp interpolation
  estimates.
\newblock {\em Comm. Math. Phys.}, 87(4):567--576, 1982/83.

\end{thebibliography}

\end{document}